\documentclass[11pt,a4paper]{article}

\usepackage{amsmath,amssymb}
\newcommand{\R}{\mathbb{R}}

\newcommand{\D}{\Delta}

\newcommand{\n}{\nabla}
\newcommand{\NNN}{\frac{N}{p}}

\newcommand{\p}{\partial}

\newcommand{\h}{\hookrightarrow}

\newtheorem{cor}{Corollary}

\newcommand{\NN}{\frac{N}{q}}
\newcommand{\N}{\frac{N}{2}}

\newcommand{\e}{\epsilon}
\newcommand{\va}{\varphi}

\newtheorem{definition}{Definition}
\newtheorem{theorem}{Theorem}
\newtheorem{notation}{Notation}
\newtheorem{proposition}{Proposition}
\newtheorem{corollaire}{Corollary}
\newtheorem{remarka}{Remark}

\newtheorem{lemme}{Lemma}

\title{Global existence of strong solution for shallow water system with large initial data on the irrotational part}
\author{Boris Haspot\thanks{Ceremade UMR CNRS 7534
Universit\'e de Paris  Dauphine,
Place du MarŽchal DeLattre De Tassigny
75775 PARIS CEDEX 16 , haspot@ceremade.dauphine.fr }}
\date{}
\begin{document}
\maketitle
\begin{abstract}
We show existence of global strong solutions  with large initial data on the irrotational part for the shallow-water system in dimension $N\geq 2$. We introduce a new notion of \textit{quasi-solutions} when the initial velocity is assumed to be irrotational, these last one exhibit regularizing effects both on the velocity and  in a very surprising way also on the density (indeed the density is a priori governed by an hyperbolic equation). We would like to point out that this smoothing effect is purely non linear and is absolutely crucial in order to deal with the pressure term as it provides new damping effects in high frequencies. In particular our result gives a first kind of answer to the problem of the existence of global weak solution for the shallow-water system. We conclude by giving new point wise decay estimates on the solution which improves the previous works \cite{HZ1,HZ2}.
\end{abstract}
\section{Introduction}
The motion of a general barotropic compressible fluid is described by the following system:
\begin{equation}
\begin{cases}
\begin{aligned}
&\p_{t}\rho+{\rm div}(\rho u)=0,\\
&\p_{t}(\rho u)+{\rm div}(\rho u\otimes u)-{\rm div}(\mu(\rho)D(u))-\n(\lambda(\rho){\rm div} u)+\n P(\rho)=0,\\
&(\rho,u)_{/t=0}=(\rho_{0},u_{0}).
\end{aligned}
\end{cases}
\label{0.1}
\end{equation}
Here $u=u(t,x)\in\R^{N}$ stands for the velocity field and $\rho=\rho(t,x)\in\R^{+}$ is the density.
The pressure $P$ is a suitable smooth function of $\rho$.
We denote by $\lambda$ and $\mu$ the two viscosity coefficients of the fluid,
which are assumed to satisfy $\mu>0$ and $\lambda+2\mu>0$. Such a condition ensures
ellipticity for the momentum equation and is satisfied in the physical cases where $\lambda+\frac{2\mu}{N}>0$. In the sequel we shall only consider the shallow-water system which corresponds to:
$$\mu(\rho)=\mu\rho\;\;\mbox{with}\;\mu>0\;\;\;\;\mbox{and}\;\;\;\;\lambda(\rho)=0.$$
We supplement the problem with initial condition $(\rho_{0},u_{0})$.
Throughout the paper, we assume that the space variable $x\in\R^{N}$ or to the periodic
box ${\cal T}^{N}_{a}$ with period $a_{i}$, in the i-th direction. We restrict ourselves to the case $N\geq2$.\\
In this paper we are interested in proving the announced result in \cite{cras}. Before giving our main result, let us recall some important results concerning the existence of strong solutions for the compressible Navier Stokes system (\ref{0.1}). The existence and uniqueness of local classical solutions for (\ref{0.1}) with smooth initial data such that the density $\rho_{0}$ is bounded
and bounded away from zero has been stated by Nash in \cite{Nash}. Let us emphasize that no stability condition was required there. On the other hand, for small smooth perturbations of a stable equilibrium with constant positive density, global well-posedness
has been proved in \cite{MN1}. More precisely  Matsumura and Nishida in \cite{MN1} obtained the existence of global strong solutions for three-dimensional polytropic ideal fluids and no outer force with initial data chosen small in the following spaces $(\rho_{0}-1,u_{0})\in H^{3}\times H^{3}$.
 Refined functional analysis has been used
during the last decades, ranging from Sobolev, Besov, Lorentz and
Triebel spaces to describe the regularity and long time behavior of
solutions to the compressible model \cite{5So}, 
, \cite{5K1}, \cite{CH}, \cite{JDE}.\\
Guided in our approach by numerous works dedicated to the incompressible Navier-Stokes equation (see e.g \cite{Meyer}):
$$
\begin{cases}
\begin{aligned}
&\p_{t}v+v\cdot\n v-\mu\D v+\n\Pi=0,\\
&{\rm div}v=0,
\end{aligned}
\end{cases}
\leqno{(NS)}
$$
we aim at solving (\ref{0.1}) in the case where the data $(\rho_{0},u_{0})$ have \textit{critical} regularity.\\
By critical, we mean that we want to solve the system (\ref{0.1}) in functional spaces with norm
 invariant by the changes of scales which leave (\ref{0.1}) invariant.
In the case of barotropic fluids, it is easy to see that the transformations:
\begin{equation}
(\rho(t,x),u(t,x))\longrightarrow (\rho(l^{2}t,lx),lu(l^{2}t,lx)),\;\;\;l\in\R,
\label{1}
\end{equation}
have that property, provided that the pressure term has been changed accordingly.\\
One of the main difficulty of compressible fluid mechanics is to deal with the vacuum, indeed in this case the momentum equation loses its parabolicity. That is why in the sequel we will work around stable equilibrium in order to control the vacuum. 
\begin{definition}
In the sequel we shall  note: $q=\rho-1$.
\end{definition}
The use of critical functional frameworks led to several new well-posedness results for compressible
fluids (see \cite{Fourier,JDE, arma, DG,CD}). In addition to have a norm invariant by (\ref{1}),
appropriate functional spaces for solving (\ref{0.1}) must provide a control on the $L^{\infty}$
norm of the density (in order to avoid vacuum and loss of parabolicity but also for dealing with the non linear term as the pressure). Additionally we also must have some Lipschitz control  on the velocity in order to control the density via the transport equation on the density. For that reason,
we restricted our study to the case where the initial data $(\rho_{0},u_{0})$ 
are in homogeneous Besov spaces such that:
$$q_{0}\in B^{\NN}_{p,1}\;\;\mbox{and}\;\;u_{0}\in B^{\frac{N}{p_{1}}-1}_{p_{1},1}$$
with $(p,p_{1})\in [1,+\infty[$ suitably chosen.\\
Recall that R. Danchin in \cite{DG} shows for the first time a result of existence of global strong solution close with \textbf{small} initial data in critical space for the scaling of the system. More precisely the initial data are chosen as follows $(q_{0},u_{0})\in (B^{\N}_{2,1}\cap B^{\N-1}_{2,1})\times B^{\N-1}_{2,1}$. The main difficulty consists in obtaining estimates on the linearized system given that the velocity and the density are coupled via the pressure. The crucial points is the obtention of damping effect on the density in order  to control the pressure term. This last result has been generalized to the case of Besov space constructed on $L^{p}$ space by F. Charve and R. Danchin in \cite{CD} and the author in \cite{arma} by using two different methods.
In \cite{CD} the authors in order to obtain estimates on the linearized system associated to (\ref{0.1}) including the convection terms combines  a paralinearization method developed by T. Hmidi and a accurate study of the linear system without the convection term in order to get estimates in Besov space. In \cite{arma} we extend the results of Charve and Danchin \cite{CD}  to the case where the Lebesgue index of Besov spaces are not the same for the density and the velocity and by reaching general index $p$ for the density. Indeed in \cite{CD} the authors need the following restriction $p<2N$ where $p$ is the Lebesgue index for the density. To do that, as in \cite{JDE} we introduce
a new notion of \textit{effective velocity}  in high frequencies which allows us to cancel out the  coupling between the velocity and the pressure. This effective velocity enables us to get as in R. Danchin in \cite{DG} a $L^{1}$ decay on $q$ in the high frequency regime. In low frequencies, the first order terms predominate, so that (\ref{0.1}) has to be treated by means
of hyperbolic energy methods (roughly $q$ and the potential part of the velocity verify a wave equation). This implies that we can treat the low regime only in spaces constructed on $L^{2}$ as it is classical that hyperbolic systems are ill-posed in general $L^{p}$ spaces. In particular the system has to be handled differently in low and high frequencies.\\
Recall that in \cite{Hoffnouv1},  Hoff  stated the existence of global weak solutions with small initial data including discontinuous initial data (namely $q_{0}\in L^{2}\cap L^{\infty}$ and is small in $L^{2}$ and $u_{0}$ is small in $L^{4}$ if $N=2$ and small
in $L^{8}$ if $N=3$). One of the major interest of the results of Hoff is to exhibit some smoothing effects on the incompressible part of the velocity $u$ and on the so-called \textit{effective viscous flux} $F=(2\mu+\lambda){\rm div}u-P(\rho)+P(\bar{\rho})$. This also plays a crucial role in the proof of  Lions for the existence of global weak solution (see \cite{1L2}). However if the results of Hoff are critical in the sense of the scaling for the density, it is not the case for the initial velocity. We would like to mention that the results of \cite{arma} makes the link in term of scaling between the results obtained in \cite{CD} and  \cite{Hoffnouv1}.\\
However the existence of global strong solution with large initial
data remains open even in dimension $N=2$ except for some very specific viscosity coefficients (see \cite{VG}). Indeed in a remarkable work  Vaigant and Kazhikhov prove the existence of global strong solution when the viscosity coefficients are chosen such that $\mu(\rho)=\mu$ and $\lambda(\rho)=\rho^{\beta}$ with $\beta>3$. To do this, they use very clever energy inequalities by taking profit of the structure of effective flux that introduce the choice of the viscosity coefficient. We would like to mention that the condition $\beta>3$ is crucial in order to get $L^{\infty}$ estimates on the density. In \cite{arma, CD,Hoffnouv1} all these works need to assume a smallness condition on the initial data in order to obtain global solution.\\
In this paper we would give a first kind of answer to this problem for a large family of initial data, more precisely we are going to prove the existence of global strong solution for initial data with large irrotational part. To do this we shall work around a irrotational \textit{quasi-solution} of the system (\ref{0.1}) (we also refer to \cite{cpde} for a such type of solution in the context of Korteweg system).
It is not clear as for the Euler system how to write the Saint-Venant system under a purely irrotational form (it means with a solution such that $u=\n\theta$) because the strong coupling between the velocity and the density, however we can check that $(\rho^{1},-\mu\n\ln\rho^{1})$ is a particular irrotational solution of the system (\ref{0.1}) (more exactly when we consider the eulerian form of (\ref{0.1})) when the pressure term is null ($P(\rho)=0$) and with:
$$
\begin{cases}
\begin{aligned}
&\p_{t}\rho^{1}-\mu\D\rho^{1}=0.\\
&\rho^{1}(0,\cdot)=\rho^{1}_{0}.
\end{aligned}
\end{cases}
$$ 
Here we recall that we assume no vacuum on the density $\rho$. It is crucial to show that $(\rho^{1},-\mu\n\ln\rho^{1})$ is well a solution of (\ref{0.1}) with $P=0$. We shall call this solution a \textit{quasi-solution} of the system (\ref{0.1}) . We recall that the control of the vacuum will be also important in order to take advantage of the parabolicity of the momentum equation.\\
It is then natural to work around this particular \textit{quasi-solution} in order to obtain global strong solution with large initial data for (\ref{0.1}), the difficulty shall consist in canceling out the effects of the pressure as the solution $(\rho^{1},-\mu\n\ln\rho^{1})$ does not take into account the pressure $P(\rho)$. One of the main argument will be to explain why by the regularizing effects on the density $\rho^{1}$ or scaling arguments this pressure can be ignored in some sense. Let us mention that the regularizing effects on $\rho^{1}$ shall play a crucial role in the proof, we shall discuss on this point in more details in the sequel.\\
We now search solution of the form $\ln\rho=\ln\rho^{1}+h^{2}$ with $\rho=\rho^{1}e^{h^{2}}$ and $u=-\mu\n\ln\rho_{1}+u^{2}$, assuming that there is no vacuum, we can rewrite the system (\ref{0.1}) under the following form:
\begin{equation}
\begin{cases}
\begin{aligned}
&\p_{t}\ln\rho+u\cdot\n\ln\rho+{\rm div}u=0,\\
&\p_{t}u+u\cdot\n -\mu\D u-\mu\n\ln\rho Du+\n F(\rho)=0,\\
&(\ln\rho,u)_{/t=0}=(\ln\rho_{0},u_{0}).
\end{aligned}
\end{cases}
\label{3systeme}
\end{equation}
By using the fact that $(\rho^{1},u^{1})=(\rho^{1},-\mu\n\ln\rho^{1})$ with:
$$
\begin{cases}
\begin{aligned}
&\p_{t}\rho^{1}-\mu\D\rho^{1}=0.\\
&\rho^{1}(0,\cdot)=\rho^{1}_{0}.
\end{aligned}
\end{cases}
$$ 
we can rewrite the system (\ref{3systeme}) as follows:
\begin{equation}
\begin{cases}
\begin{aligned}
&\p_{t}h^{2}+u\cdot\n h^{2}+{\rm div}u^{2}=-u^{2}\cdot\n \ln\rho^{1},\\
&\p_{t}u^{2}+u\cdot\n u^{2}-\mu\D u^{2}+a\n h^{2}=-a\n \ln\rho^{1}-u_{2}\cdot\n u^{1}+\mu\n\ln\rho^{1}\cdot D u^{2}\\
&\hspace{8cm}+\mu\n h^{2}\cdot D u^{1}+\mu\n h^{2}\cdot D u^{2},\\
&(q^{2},u^{2})_{/t=0}=(q^{2}_{0},u^{2}_{0}).
\end{aligned}
\end{cases}
\label{0.2}
\end{equation}
where we assume that $P(\rho)=a\rho$. In fact it would be very easy to generalize our results to general pressure, we are just interested in simplifying the notation.\\
We now are going to use the same strategy than in \cite{arma} in order to prove the existence of global strong solution $(h^{2},u^{2})$ of the system (\ref{0.2}), in particular we shall use the notion of \textit{effective velocity} introduced in \cite{JDE} and developed also in \cite{arma}, in order to obtain some damping effects on $q^{2}$. One of the main point consist in distinguishing the behavior between the low and the high frequencies as in \cite{arma}. We are going the following definition of the Besov space.
\begin{definition}
 Furthermore we will note $\widetilde{B}^{s_{1},s_{2}}_{(p_{1},r_{1}),(p_{2},r_{2})}$ the Besov space where the behavior is $B^{s_{1}}_{p_{1},r_{1}}$ in low frequencies and $B^{s_{2}}_{p_{2},r_{2}}$  in high frequencies. If $r_{1}=r_{2}$ we will simplify the notation, and we will write
 $\widetilde{B}^{s_{1},s_{2}}_{p_{1},p_{2},1}$ for $\widetilde{B}^{s_{1},s_{2}}_{(p_{1},1),(p_{2},1)}$ . For more details on the definition of these spaces we refer to the definition \ref{def1.9}.
 \end{definition}
One can now state our main result.
\begin{theorem}
Let $1\leq p<\max(4,N)$ such that $\sup(\frac{1}{q},\frac{1}{p})\leq 1+\inf(\frac{1}{p},\frac{1}{q})$, $\sup(\frac{1}{q_{1}},\frac{1}{2})\leq 1+\inf(\frac{1}{q_{1}},\frac{1}{2})$
 and  $\rho_{0}=\rho^{1}_{0}e^{h^{2}_{0}}$ and $u_{0}=-\mu \n\ln\rho^{1}_{0}+u^{2}_{0}$. Furthermore we assume that $\rho^{1}_{0}\geq c>0$, $q^{1}\in B^{0}_{\infty,1}\cap \widetilde{B}^{\frac{N}{q_{1}}-1,\NN}_{q_{1,q,\infty}}\cap \widetilde{B}^{\N-2,\NNN-2}_{2,p,1}$, $h^{2}\in \widetilde{B}^{\N-1,\NNN}_{2,p,1}$ and $u_{2}\in \widetilde{B}^{\N-1,\NNN-1}_{2,p,1}$. Then it exists  $\e$ such that if:
$$\|q^{1}_{0}\|_{\widetilde{B}^{\N-2,\NNN-2}_{2,p,1}}+\|h^{2}_{0}\|_{\widetilde{B}^{\N-1,\NNN}_{2,p,1}}+\|u^{2}_{0}\|_{\widetilde{B}^{\N-1,\NNN-1}_{2,p,1}}\leq\e,$$
then if $\frac{1}{2}+\frac{1}{p}>\frac{1}{N}$, $\frac{1}{2}+\frac{1}{q_{1}}>\frac{1}{N} $, $\frac{1}{p}+\frac{1}{q}>\frac{1}{N}$ there exists a global solution $(\rho,u)$ of the system (\ref{0.1}) written under the following form:
$\rho=\rho^{1}e^{h^{2}}\;\;\mbox{and}\;\;u=-\mu\n\ln\rho^{1}+u^{2}$ with:
\begin{equation}
\begin{cases}
\begin{aligned}
&\p_{t}\rho^{1}-\mu\D\rho^{1}=0,\\
&\rho^{1}_{t=0}=\rho^{1}_{0}.
\end{aligned}
\end{cases}
\label{achal}
\end{equation}
and such that:
$$
\begin{aligned}
&h^{2}\in \widetilde{C}(\R^{+},\widetilde{B}^{\N-1,\NNN}_{2,p,1}
)\cap L^{1}(\R^{+},
\widetilde{B}^{\N+1,\NNN}_{2,p,1})\;\;\\
&\hspace{3cm}\mbox{and}\;\;u^{2}\in \widetilde{C}(\R^{+};\widetilde{B}^{\N-1,\NNN-1}_{2,p,1})
\cap L^{1}(\R{+},\widetilde{B}^{\N+1,\NNN+1}_{2,p,1}).
\end{aligned}
$$
We refer to the definition for the definition of the hybrid Besov spaces.\label{theo1}
\end{theorem}
\begin{remarka}
We would like to emphasize on the fact that the density consists in the product of a regular function $\rho^{1}$ and of a small perturbation of the equilibrium $e^{h^{2}}$
. This point is very surprising in the sense that the density is governed by a hyperbolic equation which means that a priori we do not wait for any regularizing effects on the density. It seems that there is a singular behavior around the \textit{quasi-solution} $(\rho^{1},-\mu\n\ln\rho^{1})$ and we note that this effect is strictly non-linear. Indeed it depends on the convection term $u\cdot\n u$.\\
Furthermore this regularizing effect is crucial in order to deal with the term $u^{2}\cdot\n\ln\rho^{1}$, indeed we lose one derivate on $\n\ln\rho^{1}$, it is then crucial to obtain regularizing effects in order to treat this term.
\end{remarka}
\begin{remarka}
Up our knowledge, it is the first result of global strong solution for compressible Navier-Stokes equations with large initial data in dimension $N\geq 2$ in critical Besov spaces for the scaling of the equations. In particular it improves strongly the results of \cite{CD,arma}. Indeed the only condition of smallness on $\rho$ is on the following Besov space $\widetilde{B}^{\N-2,\NNN-2}_{2,p,1}$ which is largely subcritical, it is the same for the velocity for the Besov space $\widetilde{B}^{\N-3,\NNN-3}_{2,p,1}$. We observe that the density is not necessary small in norm $L^{\infty}$ as it is the case in \cite{arma,CD}.\\
 In particular we emphasize on the fact that the we can choose initial data with large initial data in the energy space. Indeed we recall that the energy space is such that $\rho_{0}\in L^{2}_{1}(\R^{N})$ (we refer to \cite{1L2} for the definition of the Orlicz space) and $\sqrt{\rho_{0}}u_{0}\in L^{2}(\R^{N})$. In our case in dimension $N=3$ if $p=q=2$ we need in high frequencies an assumption of smallness on $q^{1}$ in $B^{-\frac{1}{2}}_{2,1}$. And in high frequencies $L^{2}$ is embedded in $B^{-\frac{1}{2}}_{2,1}$. We can then choose initial density with large norm. It is exactly the same for the initial velocity as the smallness condition in high frequencies correspond to the norm $B^{-\frac{3}{2}}_{2,1}$.\\
In particular it gives a first kind of answer to the problem of the existence of global weak solution for the shallow-water system. Indeed in \cite{MV}, Mellet and Vasseur prove the stability of the global weak solution for the shallow-water system. However it seems extremely difficult to construct a sequence $(\rho_{n},u_{n})_{n\in\mathbb{N}}$ of regular approximate global solution of the system (\ref{0.1}) which verifies uniformly in $n$ the different entropies that use Mellet and Vasseur for proving their result. The theorem \ref{theo1} allows us to obtain the existence of global weak solution for a family of large initial data in the energy space (it means some initial density such that $\rho^{0}$ is small in high frequencies in $B^{-\frac{1}{2}}_{2,1}$ but not necessary in $L^{2}$). To do this it suffices to regularize the energy initial data to construct approximate global solution $(\rho_{n},u_{n})_{n\in\mathbb{N}}$ of the system (\ref{0.1}) and to pass to the limit by the compactness arguments developed in \cite{MV}.
\end{remarka}
\begin{remarka}
By using the same idea than in \cite{arma}, it would be possible to cancel out the coupling between $h^{2}$ and $u^{2}$, and in particular obtaining the existence of global strong solution with $h^{2}_{0}\in B^{\NNN}_{p,1}$ with $p$ close from the infinity if we assume additionaly that $\rho^{0}$ is in $H^{1}(\R^{N})$ with is conform with the energy space used in \cite{MV}. \\
Let also mention that it would be possible to improve the regularity hypothesis in low frequencies, indeed by studying only the variable $u=u^{1}+u^{2}$ in high frequencies and by working directly with $u$ in low frequencies following the same idea than in \cite{arma}.
\end{remarka}
\begin{remarka}
We would like to mention that this result is strongly related to the structure of the viscosity coefficient as we are able to construct \textit{quasi-solutions}. Indeed with constant viscosity coefficient it seems not clear how to construct \textit{quasi-solution}.
\end{remarka}
\begin{remarka}
We would like to point out that our initial density is continuous as we assume that $q^{1}_{0}$ is in $B^{0}_{\infty,1}$. In \cite{kor}, we are able to prove the existence of global strong solution with discontinuous initial density for the Korteweg system, this is due to the fact that for the Korteweg system we have also regularizing effects on $h^{2}$. Here the term $\n h^{2}.Du^{1}$ imposed a control $L^{\infty}$ on $D u^{1}$ (it means a Lipschitz control). To have this, by using the proposition \ref{asymptotique} we see that the condition $q^{1}_{0}\in B^{0}_{\infty,1}$ is optimal.
\end{remarka}
\begin{remarka}
We could weaken the condition on $(h^{2},u^{2})$ by following also the idea of \cite{CH}.
\end{remarka}
We are now going a result on the time  asymptotic behavior of our solution $(\rho,u)$.
\begin{cor}
Let $(\rho,u)$ the global strong solution constructed in the theorem \ref{theo1} with the following additional condition on the initial data, $q^{1}_{0}$ belongs to $L^{1}$. We set:
$$\e=\|h^{2}\|_{ \widetilde{L}^{\infty}(\widetilde{B}^{\N-1,\NNN}_{2,p,1})}+\|u^{2}\|_{\widetilde{L}^{\infty}(\widetilde{B}^{\N-1,\NNN-1}_{2,p,1})}.$$
Then the global solution $(\rho,u)$ satisfies the following decay estimates for:
\begin{equation}
\begin{aligned}
\|(\rho-1)(t,\cdot)\|_{L^{\infty}}
&\leq C\big(\frac{\|q^{1}_{0}\|_{L^{1}}}{(1+t)^{\N}}(1+\e)+\e\big) ,
\end{aligned}
\label{cestim5}
\end{equation}
and
\begin{equation}
\begin{aligned}
\|u(t,\cdot)\|_{B^{-1}_{\infty,\infty}}&\leq C\big(\frac{\|q^{1}_{0}\|_{L^{1}}}{(1+t)^{\N+\frac{1}{2}}}+\e),
\end{aligned}
\label{cestim5}
\end{equation}
where $C$ depends on the viscosity $\mu$ and the initial data $h^{0}_{2}$, $u_{0}^{2}$.
\label{theo2}
\end{cor}
\begin{remarka}
This result improves the works of \cite{HZ1,HZ2} as we are working with critical initial data without assuming that the regularity of $\rho_{0}$ is the same than $u_{0}$. Indeed in \cite{HZ1,HZ2} consider the compressible Navier-Stokes equations as a regularizing system of the compressible Euler system.\\
We refer also to the thesis of Rodrigues (see \cite{Ro}) for more details on the time asymptotic decay of the global strong solution with small initial data.\\
We shall come back on this problem in a forthcoming paper.
\end{remarka}
We are now going to consider the viscous shallow water model with friction. This model is also called by the french community the Saint-Venant equations and is generally used in oceanography. Indeed it allows to model vertically averaged flows  in terms of the horizontal mean velocity field $u$ and the depth variation $\rho$. In the rotating framework, the model is described by the following system::
\begin{equation}
\begin{cases}
\begin{aligned}
&\p_{t}\rho+{\rm div}(\rho u)=0,\\
&\p_{t}(\rho u)+{\rm div}(\rho u\otimes u)-{\rm div}(\mu\rho D(u))+\frac{\n \rho}{Fr^{2}}+r \rho u=0,\\
&(\rho,u)_{/t=0}=(\rho_{0},u_{0}).
\end{aligned}
\end{cases}
\label{10.1}
\end{equation}
$Fr > 0$ denotes the Froude number. The turbulent regime ($r\geq 0$) is obtained from the friction condition on the bottom, see \cite{P}. We shall work in the sequel with $r=\frac{1}{Fr^{2}}$.\\
In this case we are able to exhibit explicit solution $(\rho^{1},-\mu\n\ln\rho{1})$ which verifies (\ref{achal}), by using similar idea than in theorem \ref{theo1} we shall obtain the following result.
\begin{theorem}
Let $1\leq p<\max(4,N)$ such that $\sup(\frac{1}{q},\frac{1}{p})\leq 1+\inf(\frac{1}{p},\frac{1}{q})$, $\sup(\frac{1}{q_{1}},\frac{1}{2})\leq 1+\inf(\frac{1}{q_{1}},\frac{1}{2})$
 and  $\rho_{0}=\rho^{1}_{0}e^{h^{2}_{0}}$ and $u_{0}=-\mu \n\ln\rho^{1}_{0}+u^{2}_{0}$. Furthermore we assume that $\rho^{1}_{0}\geq c>0$, $q^{1}\in B^{0}_{\infty,1}\cap \widetilde{B}^{\frac{N}{q_{1}}-1,\NN}_{q_{1,q,\infty}}\cap \widetilde{B}^{\N-2,\NNN-2}_{2,p,1}$, $h^{2}\in \widetilde{B}^{\N-1,\NNN}_{2,p,1}$ and $u_{2}\in \widetilde{B}^{\N-1,\NNN-1}_{2,p,1}$. Then it exists  $\e$ such that if:
$$\|h^{2}_{0}\|_{\widetilde{B}^{\N-1,\NNN}_{2,p,1}}+\|u^{2}_{0}\|_{\widetilde{B}^{\N-1,\NNN-1}_{2,p,1}}\leq\e,$$
then if $\frac{1}{2}+\frac{1}{p}>\frac{1}{N}$ there exists a global solution $(\rho,u)$ of the system (\ref{0.1}) written under the following form:
$\rho=\rho^{1}e^{h^{2}}\;\;\mbox{and}\;\;u=-\mu\n\ln\rho^{1}+u^{2}$ with:
$$
\begin{cases}
\begin{aligned}
&\p_{t}\rho^{1}-\mu\D\rho^{1}=0,\\
&\rho^{1}_{t=0}=\rho^{1}_{0}.
\end{aligned}
\end{cases}
$$
and such that:
$$
\begin{aligned}
&h^{2}\in \widetilde{C}(\R^{+},\widetilde{B}^{\N-1,\NNN}_{2,p,1}
)\cap L^{1}(\R^{+},
\widetilde{B}^{\N+1,\NNN}_{2,p,1})\;\;\mbox{and}\;\;u^{2}\in \widetilde{C}(\R^{+};\widetilde{B}^{\N-1,\NNN-1}_{2,p,1})
\cap L^{1}(\R{+},\widetilde{B}^{\N+1,\NNN+1}_{2,p,1}).
\end{aligned}
$$
We refer to the definition for the definition of the hybrid Besov spaces.\label{theo3}
\end{theorem}
\begin{remarka}
Compared with the theorem \ref{theo1}, we do not need any assumption of smallness on the density $\rho^{1}_{0}$, it is completely a result of global strong solution for large initial data when $N\geq2$. It is the first result up our knowledge of global strong solution with large initial data for a compressible system.
\end{remarka}
Our paper is structured as follows. In section \ref{section2}, we give a few notation and briefly introduce the basic Fourier analysis
techniques needed to prove our result. In section \ref{section3}, we prove estimates on a linear system with convection terms. In section \ref{section4} we prove the theorems \ref{theo1}, the corollary \ref{theo2} and the theorem \ref{theo3}. Some technical continuity results for the paraproduct in hybrid Besov spaces have been postponed in appendix.
\section{Littlewood-Paley theory and Besov spaces}
\label{section2}
Throughout the paper, $C$ stands for a constant whose exact meaning depends on the context. The notation $A\lesssim B$ means
that $A\leq CB$.
For all Banach space $X$, we denote by $C([0,T],X)$ the set of continuous functions on $[0,T]$ with values in $X$.
For $p\in[1,+\infty]$, the notation $L^{p}(0,T,X)$ or $L^{p}_{T}(X)$ stands for the set of measurable functions on $(0,T)$
with values in $X$ such that $t\rightarrow\|f(t)\|_{X}$ belongs to $L^{p}(0,T)$.
Littlewood-Paley decomposition  corresponds to a dyadic
decomposition  of the space in Fourier variables.
We can use for instance any $\varphi\in C^{\infty}(\R^{N})$,
supported in
${\cal{C}}=\{\xi\in\R^{N}/\frac{3}{4}\leq|\xi|\leq\frac{8}{3}\}$
such that:
$$\sum_{l\in\mathbb{Z}}\varphi(2^{-l}\xi)=1\,\,\,\,\mbox{if}\,\,\,\,\xi\ne 0.$$
Denoting $h={\cal{F}}^{-1}\varphi$, we then define the dyadic
blocks by:
$$\D_{l}u=\varphi(2^{-l}D)u=2^{lN}\int_{\R^{N}}h(2^{l}y)u(x-y)dy\,\,\,\,\mbox{and}\,\,\,S_{l}u=\sum_{k\leq
l-1}\D_{k}u\,.$$ Formally, one can write that:
$$u=\sum_{k\in\mathbb{Z}}\D_{k}u\,.$$
This decomposition is called homogeneous Littlewood-Paley
decomposition. Let us observe that the above formal equality does
not hold in ${\cal{S}}^{'}(\R^{N})$ for two reasons:
\begin{enumerate}
\item The right hand-side does not necessarily converge in
${\cal{S}}^{'}(\R^{N})$.
\item Even if it does, the equality is not
always true in ${\cal{S}}^{'}(\R^{N})$ (consider the case of the
polynomials).
\end{enumerate}
\subsection{Homogeneous Besov spaces and first properties}
\begin{definition}
For
$s\in\R,\,\,p\in[1,+\infty],\,\,q\in[1,+\infty],\,\,\mbox{and}\,\,u\in{\cal{S}}^{'}(\R^{N})$
we set:
$$\|u\|_{B^{s}_{p,q}}=(\sum_{l\in\mathbb{Z}}(2^{ls}\|\D_{l}u\|_{L^{p}})^{q})^{\frac{1}{q}}.$$
The Besov space $B^{s}_{p,q}$ is the set of temperate distribution $u$ such that $\|u\|_{B^{s}_{p,q}}<+\infty$.
\end{definition}
\begin{remarka}The above definition is a natural generalization of the
nonhomogeneous Sobolev and H$\ddot{\mbox{o}}$lder spaces: one can show
that $B^{s}_{\infty,\infty}$ is the nonhomogeneous
H$\ddot{\mbox{o}}$lder space $C^{s}$ and that $B^{s}_{2,2}$ is
the nonhomogeneous space $H^{s}$.
\end{remarka}
\begin{proposition}
\label{derivation,interpolation}
The following properties holds:
\begin{enumerate}
\item there exists a constant universal $C$
such that:\\
$C^{-1}\|u\|_{B^{s}_{p,r}}\leq\|\n u\|_{B^{s-1}_{p,r}}\leq
C\|u\|_{B^{s}_{p,r}}.$
\item If
$p_{1}<p_{2}$ and $r_{1}\leq r_{2}$ then $B^{s}_{p_{1},r_{1}}\hookrightarrow
B^{s-N(1/p_{1}-1/p_{2})}_{p_{2},r_{2}}$.
\item $B^{s^{'}}_{p,r_{1}}\hookrightarrow B^{s}_{p,r}$ if $s^{'}> s$ or if $s=s^{'}$ and $r_{1}\leq r$.
\end{enumerate}
\label{interpolation}
\end{proposition}
Let now recall a few product laws in Besov spaces coming directly from the paradifferential calculus of J-M. Bony
(see \cite{BJM,BCD}).
\begin{proposition}
\label{produit1}
We have the following laws of product:
\begin{itemize}
\item For all $s\in\R$, $(p,r)\in[1,+\infty]^{2}$ we have:
\begin{equation}
\|uv\|_{B^{s}_{p,r}}\leq
C(\|u\|_{L^{\infty}}\|v\|_{B^{s}_{p,r}}+\|v\|_{L^{\infty}}\|u\|_{B^{s}_{p,r}})\,.
\label{2.2}
\end{equation}
\item Let $(p,p_{1},p_{2},r,\lambda_{1},\lambda_{2})\in[1,+\infty]^{2}$ such that:$\frac{1}{p}\leq\frac{1}{p_{1}}+\frac{1}{p_{2}}$,
$p_{1}\leq\lambda_{2}$, $p_{2}\leq\lambda_{1}$, $\frac{1}{p}\leq\frac{1}{p_{1}}+\frac{1}{\lambda_{1}}$ and
$\frac{1}{p}\leq\frac{1}{p_{2}}+\frac{1}{\lambda_{2}}$. We have then the following inequalities:\\
if $s_{1}+s_{2}+N\inf(0,1-\frac{1}{p_{1}}-\frac{1}{p_{2}})>0$, $s_{1}+\frac{N}{\lambda_{2}}<\frac{N}{p_{1}}$ and
$s_{2}+\frac{N}{\lambda_{1}}<\frac{N}{p_{2}}$ then:
\begin{equation}
\|uv\|_{B^{s_{1}+s_{2}-N(\frac{1}{p_{1}}+\frac{1}{p_{2}}-\frac{1}{p})}_{p,r}}\lesssim\|u\|_{B^{s_{1}}_{p_{1},r}}
\|v\|_{B^{s_{2}}_{p_{2},\infty}},
\label{2.3}
\end{equation}
when $s_{1}+\frac{N}{\lambda_{2}}=\frac{N}{p_{1}}$ (resp $s_{2}+\frac{N}{\lambda_{1}}=\frac{N}{p_{2}}$) we replace
$\|u\|_{B^{s_{1}}_{p_{1},r}}\|v\|_{B^{s_{2}}_{p_{2},\infty}}$ (resp $\|v\|_{B^{s_{2}}_{p_{2},\infty}}$) by
$\|u\|_{B^{s_{1}}_{p_{1},1}}\|v\|_{B^{s_{2}}_{p_{2},r}}$ (resp $\|v\|_{B^{s_{2}}_{p_{2},\infty}\cap L^{\infty}}$),
if $s_{1}+\frac{N}{\lambda_{2}}=\frac{N}{p_{1}}$ and $s_{2}+\frac{N}{\lambda_{1}}=\frac{N}{p_{2}}$ we take $r=1$.
\\
If $s_{1}+s_{2}=0$, $s_{1}\in(\frac{N}{\lambda_{1}}-\frac{N}{p_{2}},\frac{N}{p_{1}}-\frac{N}{\lambda_{2}}]$ and
$\frac{1}{p_{1}}+\frac{1}{p_{2}}\leq 1$ then:
\begin{equation}
\|uv\|_{B^{-N(\frac{1}{p_{1}}+\frac{1}{p_{2}}-\frac{1}{p})}_{p,\infty}}\lesssim\|u\|_{B^{s_{1}}_{p_{1},1}}
\|v\|_{B^{s_{2}}_{p_{2},\infty}}.
\label{2.4}
\end{equation}
If $|s|<\NN$ for $p\geq2$ and $-\frac{N}{p^{'}}<s<\NN$ else, we have:
\begin{equation}
\|uv\|_{B^{s}_{p,r}}\leq C\|u\|_{B^{s}_{p,r}}\|v\|_{B^{\NN}_{p,\infty}\cap L^{\infty}}.
\label{2.5}
\end{equation}
\end{itemize}
\end{proposition}
\begin{remarka}
In the sequel $p$ will be either $p_{1}$ or $p_{2}$ and in this case $\frac{1}{\lambda}=\frac{1}{p_{1}}-\frac{1}{p_{2}}$
if $p_{1}\leq p_{2}$, resp $\frac{1}{\lambda}=\frac{1}{p_{2}}-\frac{1}{p_{1}}$
if $p_{2}\leq p_{1}$.
\end{remarka}
\begin{corollaire}
\label{produit2}
Let $r\in [1,+\infty]$, $1\leq p\leq p_{1}\leq +\infty$ and $s$ such that:
\begin{itemize}
\item $s\in(-\frac{N}{p_{1}},\frac{N}{p_{1}})$ if $\frac{1}{p}+\frac{1}{p_{1}}\leq 1$,
\item $s\in(-\frac{N}{p_{1}}+N(\frac{1}{p}+\frac{1}{p_{1}}-1),\frac{N}{p_{1}})$ if $\frac{1}{p}+\frac{1}{p_{1}}> 1$,
\end{itemize}
then we have if $u\in B^{s}_{p,r}$ and $v\in B^{\frac{N}{p_{1}}}_{p_{1},\infty}\cap L^{\infty}$:
$$\|uv\|_{B^{s}_{p,r}}\leq C\|u\|_{B^{s}_{p,r}}\|v\|_{B^{\frac{N}{p_{1}}}_{p_{1},\infty}\cap L^{\infty}}.$$
\end{corollaire}
The study of non stationary PDE's requires space of type $L^{\rho}(0,T,X)$ for appropriate Banach spaces $X$. In our case, we
expect $X$ to be a Besov space, so that it is natural to localize the equation through Littlewood-Paley decomposition. But, in doing so, we obtain
bounds in spaces which are not type $L^{\rho}(0,T,X)$ (except if $r=p$).
We are now going to
define some useful spaces in which we will work, which are
a refinement of the spaces
$L_{T}^{\rho}(B^{s}_{p,r})$.
$\hspace{15cm}$
\begin{definition}
Let $\rho\in[1,+\infty]$, $T\in[1,+\infty]$ and $s_{1}\in\R$. We set:
$$\|u\|_{\widetilde{L}^{\rho}_{T}(B^{s_{1}}_{p,r})}=
\big(\sum_{l\in\mathbb{Z}}2^{lrs_{1}}\|\D_{l}u(t)\|_{L^{\rho}(L^{p})}^{r}\big)^{\frac{1}{r}}\,.$$
We then define the space $\widetilde{L}^{\rho}_{T}(B^{s_{1}}_{p,r})$ as the set of temperate distribution $u$ over
$(0,T)\times\R^{N}$ such that 
$\|u\|_{\widetilde{L}^{\rho}_{T}(B^{s_{1}}_{p,r})}<+\infty$.
\end{definition}
We set $\widetilde{C}_{T}(\widetilde{B}^{s_{1}}_{p,r})=\widetilde{L}^{\infty}_{T}(\widetilde{B}^{s_{1}}_{p,r})\cap
{\cal C}([0,T],B^{s_{1}}_{p,r})$.
Let us emphasize that, according to Minkowski inequality, we have:
$$\|u\|_{\widetilde{L}^{\rho}_{T}(B^{s_{1}}_{p,r})}\leq\|u\|_{L^{\rho}_{T}(B^{s_{1}}_{p,r})}\;\;\mbox{if}\;\;r\geq\rho
,\;\;\;\|u\|_{\widetilde{L}^{\rho}_{T}(B^{s_{1}}_{p,r})}\geq\|u\|_{L^{\rho}_{T}(B^{s_{1}}_{p,r})}\;\;\mbox{if}\;\;r\leq\rho
.$$
\begin{remarka}
It is easy to generalize proposition \ref{produit1},
to $\widetilde{L}^{\rho}_{T}(B^{s_{1}}_{p,r})$ spaces. The indices $s_{1}$, $p$, $r$
behave just as in the stationary case whereas the time exponent $\rho$ behaves according to H\"older inequality.
\end{remarka}
In the sequel we will need of composition lemma in $\widetilde{L}^{\rho}_{T}(B^{s}_{p,r})$ spaces.
\begin{lemme}
\label{composition}
Let $s>0$, $(p,r)\in[1,+\infty]$ and $u\in \widetilde{L}^{\rho}_{T}(B^{s}_{p,r})\cap L^{\infty}_{T}(L^{\infty})$.
\begin{enumerate}
 \item Let $F\in W_{loc}^{[s]+2,\infty}(\R^{N})$ such that $F(0)=0$. Then $F(u)\in \widetilde{L}^{\rho}_{T}(B^{s}_{p,r})$. More precisely there exists a function $C$ depending only on $s$, $p$, $r$, $N$ and $F$ such that:
$$\|F(u)\|_{\widetilde{L}^{\rho}_{T}(B^{s}_{p,r})}\leq C(\|u\|_{L^{\infty}_{T}(L^{\infty})})\|u\|_{\widetilde{L}^{\rho}_{T}(B^{s}_{p,r})}.$$
\item Let $F\in W_{loc}^{[s]+3,\infty}(\R^{N})$ such that $F(0)=0$. Then $F(u)-F^{'}(0)u\in \widetilde{L}^{\rho}_{T}(B^{s}_{p,r})$. More precisely there exists a function $C$ depending only on $s$, $p$, $r$, $N$ and $F$ such that:
$$\|F(u)-F^{'}(0)u\|_{\widetilde{L}^{\rho}_{T}(B^{s}_{p,r})}\leq C(\|u\|_{L^{\infty}_{T}(L^{\infty})})\|u\|^{2}_{\widetilde{L}^{\rho}_{T}(B^{s}_{p,r})}.$$
\end{enumerate}
\end{lemme}
Now we give some result on the behavior of the Besov spaces via some pseudodifferential operator (see \cite{BCD}).
\begin{definition}
Let $m\in\R$. A smooth function function $f:\R^{N}\rightarrow\R$ is said to be a ${\cal S}^{m}$ multiplier if for all muti-index $\alpha$, there exists a constant $C_{\alpha}$ such that:
$$\forall\xi\in\R^{N},\;\;|\p^{\alpha}f(\xi)|\leq C_{\alpha}(1+|\xi|)^{m-|\alpha|}.$$
\label{smoothf}
\end{definition}
\begin{proposition}
Let $m\in\R$ and $f$ be a ${\cal S}^{m}$ multiplier. Then for all $s\in\R$ and $1\leq p,r\leq+\infty$ the operator $f(D)$ is continuous from $B^{s}_{p,r}$ to $B^{s-m}_{p,r}$.
\label{singuliere}
\end{proposition}
Let us now give some estimates for the heat equation:
\begin{proposition}
\label{5chaleur} Let $s\in\R$, $(p,r)\in[1,+\infty]^{2}$ and
$1\leq\rho_{2}\leq\rho_{1}\leq+\infty$. Assume that $u_{0}\in B^{s}_{p,r}$ and $f\in\widetilde{L}^{\rho_{2}}_{T}
(B^{s-2+2/\rho_{2}}_{p,r})$.
Let u be a solution of:
$$
\begin{cases}
\begin{aligned}
&\p_{t}u-\mu\D u=f\\
&u_{t=0}=u_{0}\,.
\end{aligned}
\end{cases}
$$
Then there exists $C>0$ depending only on $N,\mu,\rho_{1}$ and
$\rho_{2}$ such that:
$$\|u\|_{\widetilde{L}^{\rho_{1}}_{T}(\widetilde{B}^{s+2/\rho_{1}}_{p,r})}\leq C\big(
 \|u_{0}\|_{B^{s}_{p,r}}+\mu^{\frac{1}{\rho_{2}}-1}\|f\|_{\widetilde{L}^{\rho_{2}}_{T}
 (B^{s-2+2/\rho_{2}}_{p,r})}\big)\,.$$
 If in addition $r$ is finite then $u$ belongs to $C([0,T],B^{s}_{p,r})$.
\end{proposition}
\begin{proposition} Let s be a positive real number and $(p,r) \in[1,+\infty]^{2}$. A constant
C exists which satisfies:
$$C^{-1}\|u\|_{B^{-2s}_{p,r}}\leq\big\|\;\|t^{s}e^{t\D}u\|_{L^{p}}\;\|_{L^{r}(\R^{+},\frac{dt}{t})}\leq C\|u\|_{B^{-2s}_{p,r}}.$$
\label{asymptotique}
\end{proposition}
\subsection{Hybrid Besov spaces}
The homogeneous Besov spaces fail to have nice inclusion properties: owing to the low frequencies, the embedding $B^{s}_{p,1}\hookrightarrow B^{t}_{p,1}$ does not hold for $s>t$. Still, the functions of $B^{s}_{p,1}$ are locally more regular than those of $B^{t}_{p,1}$: for any $\phi\in C^{\infty}_{0}$ and $u\in B^{s}_{p,1}$, the function $\phi u\in B^{t}_{p,1}$. This motivates the definition of Hybrid Besov spaces introduced by R. Danchin (see \cite{BCD,arma}) where the growth conditions satisfied by the dyadic blocks and the coefficient of integrability are not the same for low and high frequencies. Hybrid Besov spaces have been used in \cite{arma} to prove global well-posedness for compressible Navier-Stokes equation. We generalize here a little bit the definition by allowing for different Lebesgue norms in low and high frequencies.
\begin{definition}
 \label{def1.9}
Let  $l_{0}\in\mathbb{N}$, $s,t,\in\R$, $(r,r_{1})\in [1,+\infty]^{2}$ and $(p,q)\in[1,+\infty]$. We set:
$$\|u\|_{\widetilde{B}^{s,t}_{p,q,1}}=\sum_{l\leq l_{0}}2^{ls}\|\D_{l}u\|_{L^{p}}+\sum_{l> l_{0}}2^{lt}\|\D_{l}u\|_{L^{q}},$$
and:
$$\|u\|_{\widetilde{B}^{s,t}_{(p,r),(q,r_{1}})}=\big(\sum_{l\leq l_{0}}(2^{ls}\|\D_{l}u\|_{L^{p}})^{r}\big)^{\frac{1}{r}}+\big(\sum_{l> l_{0}}(2^{lt}\|\D_{l}u\|_{L^{q}})^{r_{1}}\big)^{\frac{1}{r_{1}}}.$$
\end{definition}
\begin{remarka}
It will be important in the sequel to chose $l_{0}$ big enough.
\end{remarka}
\begin{notation}
 We will often use the following notation:
$$u_{BF}=\sum_{l\leq l_{0}}\D_{l}u\;\;\;\mbox{and}\;\;\;u_{HF}=\sum_{l> l_{0}}\D_{l}u.$$
\end{notation}
\begin{remarka}
 We have the following properties:
\begin{itemize}
 \item We have $\widetilde{B}^{s,s}_{p,p,1}=B^{s}_{p,1}$.
\item If $s_{1}\geq s_{3}$ and $s_{2}\geq s_{4}$ then $\widetilde{B}^{s_{3},s_{2}}_{p,q,1}\h \widetilde{B}^{s_{1},s_{4}}_{p,q,1}$.
\end{itemize}
\end{remarka}
We shall also make use of hybrid Besov-spaces.The basic idea of paradifferentiel calculus is that
any product of two distributions $u$ and $v$ can be formally decomposed into:
$$uv=T_{u}v+T_{v}u+R(u,v)=T_{u}v+T^{'}_{v}u$$
where the paraproduct operator is defined by $T_{u}v=\sum_{q}S_{q-1}u\D_{q}v$, the remainder operator $R$, by
$R(u,v)=\sum_{q}\D_{q}u(\D_{q-1}v+\D_{q}v+\D_{q+1}v)$ and $T^{'}_{v}u=T_{v}u+R(u,v)$.\\
We recall here an important proposition on the paraproduct for hybrid Besov spaces (see \cite{arma}).
\begin{proposition}
\label{hybrid}
Let $p_{1},p_{2},p_{3},p_{4}\in[1,+\infty], (s_{1},s_{2},s_{3},s_{4})\in\R^{4}$ and $(p,q)\in[1,+\infty]^{2}$, we have then the following inequalities:
\begin{itemize}
\item If $\frac{1}{p}\leq\frac{1}{p_{2}}+\frac{1}{\lambda}\leq 1$, $\frac{1}{q}\leq\frac{1}{p_{4}}+\frac{1}{\lambda^{'}}\leq 1$ with $(\lambda,\lambda^{'})\in[1,+\infty]^{2}$ and $p_{1}\leq\lambda^{'}$, $p_{1}\leq\lambda$, $p_{3}\leq\lambda^{'}$ then:
\begin{equation}
\|T_{u}v\|_{\widetilde{B}^{s_{1}+s_{2}+\NN-\frac{N}{p_{1}}-\frac{N}{p_{2}}
,s_{3}+s_{4}+\frac{N}{q}-\frac{N}{p_{3}}-\frac{N}{p_{4}}}_{p,q,r}}\lesssim \|u\|_{\widetilde{B}^{s_{1},s_{3}}_{p_{1},p_{3},1}}\|v\|_{\widetilde{B}^{s_{2},s_{4}}_{p_{2},p_{4},r}},
\label{62}
\end{equation}
if $s_{1}+\frac{N}{\lambda^{'}}\leq\frac{N}{p_{1}}$, $s_{1}+\frac{N}{\lambda}\leq\frac{N}{p_{1}}$ and $s_{3}+\frac{N}{\lambda^{'}}\leq\frac{N}{p_{3}}$.
\item If 
$\frac{1}{q}\leq\frac{1}{p_{3}}+\frac{1}{p_{4}}$ and
$s_{3}+s_{4}+N\inf(0,1-\frac{1}{p_{3}}-\frac{1}{p_{4}})>0$
then
 \begin{equation}
 \sum_{l\geq 4}2^{l(s_{3}+s_{4}+\frac{N}{q}-\frac{N}{p_{3}}-\frac{N}{p_{4}})}\|\D_{l}R(u,v)\|_{L^{q}}\lesssim \|u\|_{\widetilde{B}^{s_{1},s_{3}}_{p_{1},p_{3},1}}\|v\|_{\widetilde{B}^{s_{2},s_{4}}_{p_{2},p_{4},r}}.
\label{63}
\end{equation}
\item If 
$\frac{1}{p}\leq\frac{1}{p_{3}}+\frac{1}{p_{4}}\leq 1$, $\frac{1}{p}\leq\frac{1}{p_{3}}+\frac{1}{p_{2}}\leq 1$, $\frac{1}{p}\leq\frac{1}{p_{1}}+\frac{1}{p_{4}}\leq 1$, $\frac{1}{p}\leq\frac{1}{p_{1}}+\frac{1}{p_{2}}\leq 1$ and
$s_{3}+s_{4}>0$, $s_{3}+s_{2}>0$, $s_{4}+s_{1}>0$, $s_{1}+s_{2}>0$
then
 \begin{equation}
 \sum_{l\leq 4}2^{l(s_{1}+s_{2}+\frac{N}{p}-\frac{N}{p_{1}}-\frac{N}{p_{2}})}\|\D_{l}R(u,v)\|_{L^{p}}\lesssim \|u\|_{\widetilde{B}^{s_{1},\frac{N}{p_{3}}-\frac{N}{p_{1}}+s_{1}}_{p_{1},p_{3},1}}\|v\|_{\widetilde{B}^{s_{2},\frac{N}{p_{4}}-\frac{N}{p_{2}}+s_{2}}_{p_{2},p_{4},r}}.
\label{63}
\end{equation}
with $s_{3}=\frac{N}{p_{3}}-\frac{N}{p_{1}}+s_{1}$ and $s_{4}=\frac{N}{p_{4}}-\frac{N}{p_{2}}+s_{2}$.
\item If $u\in L^{\infty}$, we also have:
\begin{equation}
\|T_{u}v\|_{\widetilde{B}^{s_{1}
,s_{2}}_{p,q,r}}\lesssim \|u\|_{L^{\infty}}\|v\|_{\widetilde{B}^{s_{1},s_{2}}_{p,q,r}},
\label{65}
\end{equation}
and if $\min(s_{1},s_{2})>0$ then:
\begin{equation}
\|R(u,v)\|_{\widetilde{B}^{s_{1}
,s_{2}}_{p,q,r}}\lesssim \|u\|_{L^{\infty}}\|v\|_{\widetilde{B}^{s_{1},s_{2}}_{p,q,r}}.
\label{66}
\end{equation}
\end{itemize}
\end{proposition}
\section{A linear model with convection}
\label{section4}
In this section, we will explain how we deal with the linear system associated to the system (\ref{0.2}) that we can write under the following form:
\begin{equation}
\begin{cases}
\begin{aligned}
&\p_{t}h^{2}+v\cdot\n h^{2}+{\rm div}u^{2}=F-u^{2}\cdot\n\ln\rho^{1},\\
&\p_{t}u^{2}+v\cdot\n u^{2}-\mu\D u^{2}+a\n h^{2}=G-u_{2}\cdot\n u^{1}+\mu\n\ln\rho^{1}\cdot D u^{2}\\
&\hspace{9cm}+\mu\n h^{2}\cdot D u^{1},\\
&(h^{2},u^{2})_{/t=0}=(h^{2}_{0},u^{2}_{0}).
\end{aligned}
\end{cases}
\label{0.2a}
\end{equation}
Here $(F,G)$ are external force with a regularity that we shall precise in the proposition (\ref{Danchinbas}) and  $(q^{1},u^{1}, v)$ are such that:
$$q^{1}\in \widetilde{L}^{\infty}(B^{\NN}_{p,1})\cap \widetilde{L}^{1}(B^{\NN+2}_{p,1}),\;u^{1},v\in \widetilde{L}^{\infty}(B^{\NN-1}_{p,1})\cap \widetilde{L}^{1}(B^{\NN+1}_{2,1}),$$
and:
$$0<c\leq \rho^{1}\leq M,$$
with $p$ verifying the conditions of theorem \ref{theo1}. We would like to start with recalling the following result which has been proved in \cite{arma,DG} by using two different method. In particular in \cite{arma}, we introduce the notion of \textit{effective velocity}.
\begin{proposition}
Let $p\leq\max(4,N)$. Let $s=\NNN$ and $s^{'}=\N-1$. Let $(\rho^{1},u^{1})=(1,0)$ and $(h^{2},u^{2})$ the solution of (\ref{0.2a}). There exists a constant $C$ depending only on $\mu$, $N$, $s$ and $s^{'}$ such that the following estimate holds:
$$
\begin{aligned}
&\|(h^{2},u^{2}(t)\|_{\widetilde{B}^{s^{'}-1,s}_{2,p,1}\times \widetilde{B}^{s^{'}-1,s-1}_{2,p,1}}+\int^{t}_{0}\|(h^{2},u^{2})(s)\|_{\widetilde{B}^{s^{'}+1,s}_{2,p,1}\times \widetilde{B}^{s^{'}+1,s+1}_{2,p,1}}ds \\
& \hspace{2cm}\leq C \big(\|(h^{2}_{0},u^{2}_{0})\|_{\widetilde{B}^{s^{'}-1,s}_{2,p,1}\times \widetilde{B}^{s^{'}-1,s-1}_{2,p,1}}
+\int^{t}_{0}e^{-V(s)}\|(F,G)(s)\|_{\widetilde{B}^{s^{'}-1,s}_{2,p,1}\times \widetilde{B}^{s^{'}-1,s-1}_{2,p,1}}ds\big).
\end{aligned}
$$
with $V(T)=\int^{T}_{0}\|\n v(s)\|_{L^{\infty}}ds$.
\label{propestim}
\end{proposition}
We are going to generalize this proposition to the case when $(q^{1},u^{1})$ is different from $(0,0)$. More precisely we have the following proposition:
\begin{proposition}
\label{Danchinbas}
Let $p\leq\max(4,N)$, $q$ such that $\frac{N}{p}-1\leq\frac{N}{q}$ and $(F,G)\in \widetilde{L}^{1}(B^{s-1,s}_{2,p,1})\times  \widetilde{L}^{1}(B^{s-1}_{2,p,1})$ and  $(h^{2},u^{2})$ a solution of (\ref{0.2a}), let $s=\NNN$ and $s^{'}=\N$. The following estimate holds:
$$
\begin{aligned}
&\|(h^{2},u^{2})\|_{\widetilde{L}_{T}^{\infty}(\widetilde{B}^{s^{'}-1,s}_{2,p,1}\times \widetilde{B}^{s^{'}-1,s-1}_{2,p,1})}+\|(h^{2},u^{2})\|_{\widetilde{L}^{1}_{T}(\widetilde{B}^{s^{'}+1,s}_{2,p,1}\times \widetilde{B}^{s^{'}+1,s+1}_{2,p,1})} \\
& \hspace{4cm}\leq C e^{V(T)}\big(\|(h^{2}_{0},u^{2}_{0})\|_{\widetilde{B}^{s^{'}-1,s}_{2,p,1}\times \widetilde{B}^{s^{'}-1,s-1}_{2,p,1}}
+\|(F,G)\|_{\widetilde{L}^{1}(\widetilde{B}^{s^{'}-1,s}_{2,p,1}\times \widetilde{B}^{s^{'}-1,s-1}_{2,p,1})}\big).
\end{aligned}
$$
with $V(T)=\int^{T}_{0}(\|q^{1}(s)\|^{4}_{B^{\frac{N}{q_{1}}-\frac{1}{2},\NN+\frac{1}{2}}_{q_{1},q,\infty}}+\|q^{1}(s)\|_{ \widetilde{B}^{\frac{N}{q_{1}}+1,\NN+2}_{q_{1},q,\infty}}+\|\n u^{1}(s)\|_{L^{\infty}\cap B^{\frac{N}{q_{1}}-1,\NN}_{q_{1},q,\infty}}+\|\n v(s)\|_{L^{\infty}})ds$.
\end{proposition}
{\bf Proof:} By using the proposition \ref{propestim} we obtain that:
\begin{equation}
\begin{aligned}
&\|(h^{2},u^{2})\|_{\widetilde{L}_{T}^{\infty}(\widetilde{B}^{s^{'}-1,s}_{2,p,1}\times \widetilde{B}^{s^{'}-1,s-1}_{2,p,1})}+\|(h^{2},u^{2})\|_{\widetilde{L}^{1}_{T}(\widetilde{B}^{s^{'}+1,s}_{2,p,1}\times \widetilde{B}^{s^{'}+1,s+1}_{2,p,1})} \\
&\hspace{0,5cm}\leq C e^{V(T)}\big(\|(h^{2}_{0},u^{2}_{0})\|_{\widetilde{B}^{s^{'}-1,s}_{2,p,1}\times \widetilde{B}^{s^{'}-1,s-1}_{2,p,1}}
+\|(F_{1},G_{1})\|_{\widetilde{L}^{1}(\widetilde{B}^{s^{'}-1,s}_{2,p,1}\times \widetilde{B}^{s^{'}-1,s-1}_{2,p,1})}\big).
\end{aligned}
\label{1estimimp}
\end{equation}
with:
$$
\begin{aligned}
&F_{1}=F-u^{2}\cdot\n\ln\rho^{1},\\
&G_{1}=G-u_{2}\cdot\n u^{1}+\mu\n\ln\rho^{1}\cdot D u^{2}+\mu\n h^{2}\cdot D u^{1}.
\end{aligned}
$$
Therefore, it is only a matter of proving appropriate estimates for $F_{1}$, $G_{1}$
 by using properties of continuity on the paraproduct (we refer in particular to the appendix of \cite{arma} when the Besov spaces are hybrid).\\
 We start with the first term $u^{2}\cdot\n \ln\rho^{1}$, we have then by proposition \ref{hybrid}, interpolation (with $s+\frac{1}{2}=\frac{1}{4}(s-1)+\frac{3}{4}(s+1)$), composition estimates and H\"older inequality:
 \begin{equation}
 \begin{aligned}
 &\int^{T}_{0} \|u^{2}\cdot\n \ln\rho^{1}(s)\|_{\widetilde{B}^{s^{'}-1,s}_{2,p,1}}ds\leq C\int^{T}_{0}(\|u^{2}(s)\|_{\widetilde{B}^{s^{'}-1,s-1}_{2,1}}\|\n \ln\rho^{1}(s)\|_{\widetilde{B}^{\frac{N}{q_{1}},\NN+1}_{q_{1},q,\infty}}\\
 &\hspace{5cm}+\|u^{2}(s)\|_{\widetilde{B}^{s^{'}+\frac{1}{2},s+\frac{1}{2}}_{2,p,1}}\|\n \ln\rho^{1}(s)\|_{\widetilde{B}^{\frac{N}{q_{1}}-\frac{3}{2},\NN-\frac{1}{2}}_{q_{1},q,\infty}})
 ds,\\
 &\leq \int^{T}_{0}(\|u^{2}(s)\|_{\widetilde{B}^{s^{'}-1,s-1}_{2,p,1}}\|q^{1}(s)\|_{B^{\frac{N}{q_{1}}+1,\NN+2}_{q_{1},q,\infty}}+\|u^{2}(s)\|_{\widetilde{B}
 ^{s^{'}+\frac{1}{2},s+\frac{1}{2}}_{2,p,1}}\|q^{1}(s)\|_{B^{\frac{N}{q_{1}}-\frac{1}{2}, \NN+\frac{1}{2}}_{q_{1},q,\infty}})
 ds,\\
 &\leq \int^{T}_{0}(\|u^{2}(s)\|_{\widetilde{B}^{s^{'}-1,s-1}_{2,p,1}}\|q^{1}(s)\|_{B^{\frac{N}{q_{1}}+1,\NN+2}_{q_{1},q,\infty}}\\
 &\hspace{3cm}+\|u^{2}(s)\|^{\frac{1}{4}}_{\widetilde{B}
 ^{s^{'}-1,s-1}_{2,p,1}}\|u^{2}(s)\|^{\frac{3}{4}}_{B^{s^{'}+1,s+1}_{2,p,1}}\|q^{1}(s)\|_{\widetilde{B}^{\frac{N}{q_{1}}-\frac{1}{2},\NN+\frac{1}{2}}_{q_{1},q,\infty}})
 ds,\\
 &\leq \int^{T}_{0}(\|u^{2}(s)\|_{\widetilde{B}^{s^{'}-1,s-1}_{2,p,1}}\|q^{1}(s)\|_{B^{\frac{N}{q_{1}}+1,\NN+2}_{q_{1},q,\infty}}\\
 &\hspace{2cm}+
 \frac{1}{4\e}\|u^{2}(s)\|_{\widetilde{B}^{s^{'}-1,s-1}_{2,p,1}}\|q^{1}(s)\|^{4}_{B^{\frac{N}{q_{1}}-\frac{1}{2},\NN+\frac{1}{2}}_{q_{1},q,\infty}}
+\frac{3\e^{\frac{1}{3}}}{4} \|u^{2}(s)\|_{\widetilde{B}^{s^{'}+1,s+1}_{2,p,1}}
 ds,
\end{aligned}
\label{1estim}
\end{equation}
We are going to deal with the terms of $G_{1}$, let us start with the term $u^{2}\cdot\n u^{1}$ , we have as $-\NN<s-1<\NN$ and as $\n u^{1}\in L^{1}_{T}(L^{\infty})$:
\begin{equation}
\int^{T}_{0} \|u^{2}\cdot\n u^{1}(s)\|_{\widetilde{B}^{s^{'}-1,s-1}_{2,p,\infty}}ds\leq  \int^{T}_{0}\|u_{2}\|_{\widetilde{B}^{s^{'}-1,s-1}_{2,p,1}}\|\n u^{1}(s)\|_{\widetilde{B}^{\frac{N}{q_{1}}-1,\NN}_{q_{1},q,\infty}\cap L^{\infty}}ds.
\label{estim2}
\end{equation}
We now have by H\"older inequality:
\begin{equation}
\begin{aligned}
&\int^{T}_{0} \|\n \ln\rho^{1}\cdot D u^{2}(s)\|_{\widetilde{B}^{s^{'}-1,s-1}_{2,p,1}}ds\leq \int^{T}_{0} \|\n\ln q^{1}\|_{\widetilde{B}^{\frac{N}{q_{1}}-\frac{3}{2},\NN-\frac{1}{2}}_{q_{1},q,\infty}} \|D u^{2}(s)\|_{\widetilde{B}^{s^{'}-\frac{1}{2},s-\frac{1}{2}}_{2,p,1}}\\
&\hspace{5cm}+\|\n\ln q^{1}\|_{\widetilde{B}^{\frac{N}{q_{1}}-\frac{1}{2},\NN+\frac{1}{2}}_{q_{1},q,\infty}} \|D u^{2}(s)\|_{\widetilde{B}^{s^{'}-\frac{3}{2},s-\frac{3}{2}}_{2,p,1}}ds\\
&\leq \int^{T}_{0} (\frac{1}{\e}\|\n q^{1}\|^{2}_{\widetilde{B}^{\frac{N}{q_{1}}-1,\NN}_{q_{1},q,1}} \|u^{2}(s)\|_{\widetilde{B}^{s^{'}-1,s-1}_{2,p,1}}+\e \|u^{2}(s)\|_{\widetilde{B}^{s^{'}+1,s+1}_{2,p,1}})
ds
\end{aligned}
\label{aestim2}
\end{equation}
It remains to deal with $\n h^{2}\cdot D u^{1}$ and by using the fact that $\frac{N}{q}$:
\begin{equation}
\int^{T}_{0} \|\n\ln h^{2}\cdot D u^{1}\|_{\widetilde{B}^{s^{'}-1,s-1}_{2,p,1}}\leq  \int^{T}_{0}\|h^{2}\|_{\widetilde{B}^{s^{'}-1,s}_{2,p,1}}\|D u^{1}(s)\|_{\widetilde{B}^{\frac{N}{q_{1}}-1,\NN}_{q_{1},q,1}}ds,
\label{estim3}
\end{equation}
By combining the estimates (\ref{1estimimp}), (\ref{1estim}), (\ref{estim2}), (\ref{aestim2}), (\ref{estim3}) and Gronwall lemma, we achieve the proposition.
\section{The proof of theorem \ref{theo1}}
\label{section5}
\subsection{Proof of the existence}
We recall here that $(\rho^{1},u^{1})=(\rho^{1},-\mu\n\ln\rho^{1})$ is a quasi-solution of the system (\ref{0.1}) with:
\begin{equation}
\begin{aligned}
&\p_{t}\rho^{1}-\mu\D\rho^{1}=0,\\
&\rho^{1}_{t=0}=\rho^{1}_{0}.
\end{aligned}
\end{equation}
It means that $(\rho^{1},u^{1})$ is an irrotational solution of the approximate system:
\begin{equation}
\begin{cases}
\begin{aligned}
&\p_{t}\rho+{\rm div}(\rho u)=0,\\
&\p_{t}(\rho u)+{\rm div}(\rho u\otimes u)-{\rm div}(\mu\rho D(u))=0,\\
&(\rho,u)_{/t=0}=(\rho^{1}_{0},u^{1}_{0}).
\end{aligned}
\end{cases}
\label{0.1bis}
\end{equation}
As we assume that $q_{0}^{1}$ is in $\widetilde{B}^{\frac{N}{q_{1}}-1,\NN}_{q_{1},q,\infty}\cap B^{0}_{\infty,1}\cap \widetilde{B}^{\N-2,\NNN-2}_{2,p,1}$, by using proposition \ref{chaleur} and the fact that $q^{1}$ verifies an heat equation, we show that for any $T>0$:
\begin{equation}
\begin{aligned}
&\|q^{1}\|_{\widetilde{L}^{\infty}(\widetilde{B}^{\frac{N}{q_{1}}-1,\NN}_{q_{1},q,\infty}\cap B^{0}_{\infty,1})}+\|q^{1}\|_{\widetilde{L}^{1}(\widetilde{B}^{\frac{N}{q_{1}}+1,\NN+2}_{q_{1},q,\infty}\cap B^{2}_{\infty,1})}\leq
C\|q^{1}_{0}\|_{\widetilde{B}^{\frac{N}{q_{1}}-1,\NN}_{q_{1},q,\infty}\cap B^{0}_{\infty,1}},\\
&\|u^{1}\|_{\widetilde{L}^{\infty}(\widetilde{B}^{\frac{N}{q_{1}}-2,\NN-1}_{q_{1},q,\infty}\cap B^{-1}_{\infty,1})}+\|q^{1}\|_{\widetilde{L}^{\infty}(\widetilde{B}^{\frac{N}{q_{1}},\NN+1}_{q_{1},q,\infty}\cap B^{1}_{\infty,1})}\leq
C\|u^{1}_{0}\|_{\widetilde{B}^{\frac{N}{q_{1}}-1,\NN}_{q_{1},q,\infty}\cap B^{0}_{\infty,1}},
\end{aligned}
\label{chal1}
\end{equation}
and:
\begin{equation}
\begin{aligned}
&\|q^{1}\|_{\widetilde{L}^{\infty}(\widetilde{B}^{\frac{N}{2}-2,\NNN-2}_{2,p,1})}+\|q^{1}\|_{\widetilde{L}^{1}(\widetilde{B}^{\frac{N}{2},\NNN}_{2,p,1})}\leq
C\|q^{1}_{0}\|_{\widetilde{B}^{\frac{N}{2}-2,\NNN-2}_{2,p,1}},\\
\end{aligned}
\label{chal2}
\end{equation}
By the maximum principle, we also obtain that for all $(t,x)\in\R^{+}\times\R^{N}$:
\begin{equation}
0\leq c\leq\rho^{1}(t,x)\leq \|\rho^{1}_{0}\|_{L^{\infty}}
\label{vide}
\end{equation}
Our goal now consists in solving the system (\ref{0.2}) in order to obtain solution of the system (\ref{0.1}) under the form
$\rho=\rho^{1}\e^{h^{2}}$, $u=u^{1}+u^{2}$ with $\rho\geq c^{'}>0$. To do this, we use a standard scheme:
\begin{enumerate}
\item We smooth out the data and get a sequence of local solutions $(h^{2}_{n},u^{2}_{n})_{n\in\mathbb{N}}$ on $[0,T_{n}]$ to (\ref{0.2})
 by using the result of \cite{JDE}.
\item 
We prove uniform estimates on $(h^{2}_{n},u^{2}_{n})$ on $[0,T_{n}]$ by using the proposition \ref{Danchinbas} and we deduce that $T_{n}=+\infty$ . 
\item We use compactness to prove that the sequence $(h^{2}_{n},u^{2}_{n})$ converges, up to extraction, to a solution of (\ref{0.2}).
\end{enumerate}
\subsubsection*{Construction of approximate solutions}
We smooth out the data as follows:
$$(h_{0}^{2})_{n}=S_{n}h_{0}^{2}\;\;\mbox{and}\;\;(u_{0}^{2})_{n}=S_{n}u_{0}^{2}\;\;\;.$$
Note that we have:
$$\forall l\in\mathbb{Z},\;\;\|\D_{l}(h_{0}^{2})_{n}\|_{L^{p}}\leq\|\D_{l}h_{0}^{2}\|_{L^{p}}\;\;\;\mbox{and}\;\;\;\|(h_{0}^{2})_{n}\|
_{\widetilde{B}^{\N-1,\NNN}_{2,p,1}}\leq \|h_{0}^{2}\|_{\widetilde{B}^{\N-1,\NNN}_{2,p,1}},$$
and similar properties for $(u_{0}^{2})_{n}$, a fact which will be used repeatedly during the next
steps. Now, according \cite{JDE}, one can solve (\ref{0.2}) with the smooth data $((q_{0}^{2})_{n},(u_{0}^{2})_{n})$.
We get a solution $(h^{2}_{n},u^{2}_{n})$ on a non trivial time interval $[0,T_{n}]$
such that:
\begin{equation}
\begin{aligned}
&h^{2}_{n}\in\widetilde{C}([0,T_{n}],B^{N}_{2,1}\cap B^{\N-1}_{2,1})\;\;u^{2}_{n}\in\widetilde{C}([0,T_{n}],,B^{\N-1}_{2,1})\cap
\widetilde{L}^{1}
([0,T_{n}],,B^{\N+1}_{2,1}).
\end{aligned}
\label{a26}
\end{equation}
\subsubsection*{Uniform bounds}
In the sequel we set:
$$\ln\rho_{n}=\ln\rho^{1}+h^{2}_{n}\;\;\mbox{and}\;\;u_{n}=u^{1}+u^{2}_{n}.$$
We recall that $(h^{2}_{n},u^{2}_{n})$ satisfies the following system:
\begin{equation}
\begin{cases}
\begin{aligned}
&\p_{t}h^{2}_{n}+u_{n}\cdot\n h^{2}_{n}+{\rm div}u^{2}_{n}=u^{2}_{n}\cdot\n\ln\rho^{1},\\
&\p_{t}u^{2}_{n}+u_{n}\cdot\n u_{n}^{2}-\mu\D u_{n}^{2}+a\n h_{n}^{2}=-u^{2}_{n}\cdot\n u^{1}+\mu\n\ln\rho^{1}\cdot D u_{n}^{2}\\
&\hspace{5cm}-a\n\ln\rho^{1}+\mu\n h_{n}^{2}\cdot D u^{1}+\mu\n h_{n}^{2}\cdot D u_{n}^{2},\\
&(h^{2}_{n},u_{n}^{2})_{/t=0}=((h^{2}_{0})_{n},(u^{2}_{0})_{n}).
\end{aligned}
\end{cases}
\label{0.2aa}
\end{equation}
In this part, we aim at getting uniform estimates on $(h^{2}_{n},u_{n}^{2})$ in the following space $F_{T}$ with the norm $\|\cdot\|_{F_{T}}$:
$$
\begin{aligned}
F_{T}=\big(\widetilde{L}_{T}^{\infty}(\widetilde{B}^{\N-1,\NNN}_{2,p,1})\cap \widetilde{L}_{T}^{1}(\widetilde{B}^{\N+1,\NNN}_{2,p,1})\big)\times
\big(\widetilde{L}^{\infty}_{T}(\widetilde{B}^{\N-1,\frac{N}{p}-1}_{2,p,1})
\cap \widetilde{L}_{T}^{1}(\widetilde{B}^{\N+1,\NNN+}_{2,p,1}\big).
\end{aligned}
$$
$$
\begin{aligned}
\|(h^{2}_{n},u_{n}^{2})\|_{F_{T}}=\|h^{2}_{n}\|_{\widetilde{L}_{T}^{\infty}(\widetilde{B}^{\N-1,\NNN}_{2,p,1})\cap \widetilde{L}_{T}^{1}(\widetilde{B}^{\N+1,\NNN}_{2,p,1})}+\|u_{n}^{2}\|_{
\widetilde{L}^{\infty}_{T}(\widetilde{B}^{\N-1,\frac{N}{p}-1}_{2,p,1})
\cap \widetilde{L}_{T}^{1}(\widetilde{B}^{\N+1,\NNN+1}_{2,p,1})}.
\end{aligned}
$$
We can observe that $(h^{2}_{n},u^{2}_{n})$ verifies exactly the system (\ref{0.2a}) with $v=u^{n}$ and:
$$
\begin{aligned}
&F^{n}=0,\\
&G^{n}=-a\n\ln\rho^{1}+\mu\n h_{n}^{2}\cdot D u^{1}+\mu\n h_{n}^{2}\cdot D u_{n}^{2}.
\end{aligned}
$$
By using proposition \ref{Danchinbas}, we have the following estimate on $(h^{2}_{n},u^{2}_{n})$:
$$
\begin{aligned}
&\|(h_{n}^{2},u^{2}_{n})\|_{F_{T}}\leq C e^{V^{n}_{1}(T)}\big(\|(h^{2}_{0})_{n}\|_{\widetilde{B}^{\N-1,\NNN}_{2,p,1}}+\|(u^{2}_{0})_{n}\|_{\widetilde{B}^{\N-1,\frac{N}{p}-1}_{2,p,1}}
+\|G^{n}\|_{\widetilde{L}_{T}^{1}(\widetilde{B}^{\N-1,\NNN-1}_{2,p,1})}\big).
\end{aligned}
$$
Therefore, it is only a matter of proving appropriate estimates for $G_{1}^{n}$ by using proposition \ref{hybrid}.
We begin by estimating $h_{n}^{2}\cdot D u_{n}^{2}$ and we obtain:
\begin{equation}
\|\n h_{n}^{2}\cdot D u_{n}^{2}\|_{\widetilde{L}_{T}^{1}(\widetilde{B}^{\N-1,\NNN-1}_{2,p,1})}\leq 
\|\n h_{n}^{2}\|_{\widetilde{L}_{T}^{\infty}(\widetilde{B}^{\N-1,\NNN-1}_{2,p,1})}\| D u_{n}^{2}\|_{\widetilde{L}_{T}^{1}(\widetilde{B}^{\N,\NNN}_{2,p,1})}.
\label{estimb1}
\end{equation}
Similarly, we obtain as $\N1\leq\frac{N}{q_{1}}$ and  $\NNN-1\leq\frac{N}{q}$:
\begin{equation}
\begin{aligned}
&\|\n h_{n}^{2}\cdot D u^{1}\|_{\widetilde{L}_{T}^{1}(\widetilde{B}^{\N-1,\NNN-1}_{2,p,1})}\leq 
\|\n h_{n}^{2}\|_{\widetilde{L}_{T}^{\infty}(\widetilde{B}^{\N-1,\NNN-1}_{2,p,1})}\| D u^{1}\|_{\widetilde{L}_{T}^{1}(\widetilde{B}^{\frac{N}{q_{1}},\NN}_{2,q,1})}
\end{aligned}
\label{estimb2}
\end{equation}
And finally by using (\ref{chal2}), we have:
\begin{equation}
\|q^{1}\|_{\widetilde{L}^{\infty}(\widetilde{B}^{\frac{N}{2}-2,\NNN-2}_{2,p,1})}+\|q^{1}\|_{\widetilde{L}^{1}(\widetilde{B}^{\frac{N}{2},\NNN}_{2,p,1})}\leq
C\|q^{1}_{0}\|_{\widetilde{B}^{\frac{N}{2}-2,\NNN-2}_{2,p,1}},
\label{estimb3}
\end{equation}
We would like to point out  the fact that $\|q^{1}_{0}\|_{\widetilde{B}^{\frac{N}{2}-2,\NNN-2}_{2,p,1}}$ small is mandatory in order to deal with the bootstrap argument. Here we use strongly the regularizing effect on $q^{1}$.\\
From a standard bootstrap argument, it is now easy to conclude that there exists a positive constant $c$ such that if the data has been chosen so small as to satisfy:
$$\|q^{1}_{0}\|_{\widetilde{B}^{\N-2,\NNN-2}_{2,p,1}}+\|h^{2}_{0}\|_{\widetilde{B}^{\N-1,\NNN}_{2,p,1}}+\|u^{2}_{0}\|_{\widetilde{B}^{\N-1,\NNN-1}_{2,p,1}}\leq\e,$$
then $T_{n}=+\infty$. Furthermore it exists $C>0$ such that for all $t\in \R$:
$$\|(q^{n},u^{n})\|_{F_{t}}\leq C,\;\;\forall t\in\R.$$
\subsubsection*{Compactness arguments}
Let us first focus on the convergence of $(h_{n}^{2})_{n\in\mathbb{N}}$. We claim that, up to extraction, $(h^{2}_{n})_{n\in\mathbb{N}}$ converges in the distributional sense to some function $h$ such that:
\begin{equation}
 h\in\widetilde{L}^{\infty}(\widetilde{B}^{\N-1,\NNN}_{2,p,1})\cap \widetilde{L}^{1}(\widetilde{B}^{\N+1,\NNN}_{2,p,1}).
\label{49}
\end{equation}
The proof is based on Ascoli's theorem and compact embedding for Besov spaces. As similar arguments have been employed in \cite{BCD}, we only give the outlines of the proof.
We may write that:
$$\p_{t}h^{2}_{n}+u_{n}\cdot\n h^{2}_{n}+{\rm div}u^{2}_{n}=u^{2}_{n}\cdot\n\ln\rho^{1},
.$$
Since $(u^{2}_{n})_{n\in\mathbb{N}}$ is uniformly bounded in $\widetilde{L}^{2}(\widetilde{B}^{\N,\NNN}_{2,p,1})$
and $\n\ln\rho^{1}\in \widetilde{L}^{\infty}(\widetilde{B}^{\frac{N}{q_{1}}-2,\NN-1}_{2,q,1})$, we have $u^{2}_{n}\cdot\n\ln\rho^{1}$ which is bounded in
$\widetilde{L}^{2}(\widetilde{B}^{\N-1,\frac{N}{p}-1}_{2,p,1})$. Similarly $u^{n}\cdot\n h_{n}^{2}$ is bounded in $\widetilde{L}^{2}(\widetilde{B}^{\N-1,\frac{N}{p}-1}_{2,p,1})$. Finally  we have proved that $\p_{t}h^{2}_{n}$ is bounded in $\widetilde{L}^{2}(\widetilde{B}^{\N-1,\frac{N}{p}-1}_{2,p,1})$, it means that $(h^{2}_{n})_{n\in\mathbb{N}}$  seen as a sequence of $\widetilde{B}^{\N-1,\frac{N}{p}-1}_{2,p,1}$ valued functions is locally equicontinuous in $\R^{+}$.
In addition $(h^{2}_{n})_{n\in\mathbb{N}}$ is bounded in $C(\R^{+},\widetilde{B}^{\N-1,\frac{N}{p}-1}_{2,p,1})$. As the embedding $\widetilde{B}^{\N-1,\frac{N}{p}-1}_{2,p,1}\cap \widetilde{B}^{\N,\frac{N}{p}}_{2,p,1}$ is locally compact (see \cite{BCD}, Chap2), one can thus conclude by means of Ascoli's theorem and Cantor diagonal extraction process that there exists some distribution $h$ such that up to an omitted extraction $(\psi h^{2}_{n})_{n\in\mathbb{N}}$ converges to $\psi h$
in $C(\R^{+},\widetilde{B}^{\N-1,\frac{N}{p}-1}_{2,p,1})$ for all smooth $\psi$ with compact support in $\R^{+}\times\R^{N}$. Then by using the so-called Fatou property for the Besov spaces, one can conclude that (\ref{49}) is satisfied. (the reader may consult \cite{BCD}, Chap 10 too).
By proceeding similarly, we can prove that up to extraction, $(u^{2}_{n})_{n\in\mathbb{N}}$ converges in the distributional sense to some function $u^{2}$ such that:
\begin{equation}
 u^{2}\in\widetilde{L}^{\infty}(\widetilde{B}^{\N-1,\frac{N}{p}-1}_{2,p_{1},1})\cap \widetilde{L}^{1}(\widetilde{B}^{\N+1,\NNN+1}_{2,p,1}).
\label{51}
\end{equation}
In order to complete the proof of the existence part of theorem \ref{theo1}, it is only a matter of checking the continuity properties with respect to time, namely that:
$$
\begin{aligned}
&h^{2}\in\widetilde{C}(\R^{+},\widetilde{B}^{\N-1,\NNN}_{2,p,1})\;\;\mbox{and}\;\;u^{2}\in \widetilde{C}(\R^{+},\widetilde{B}^{\N-1,\NNN}_{2,p,1}).
\end{aligned}
$$
As regards $h^{2}$, it suffices to notice that, according to (\ref{49}), (\ref{51}) and to the product laws in the Besov spaces, we have:
$$\p_{t}h^{2}+u\cdot\n h^{2}+{\rm div}u^{2}=u^{2}\cdot\n\ln\rho^{1}\in \widetilde{L}^{1}(\widetilde{B}^{\N,\NNN}_{2,p,1}).$$
As $h^{2}_{0}\in \widetilde{B}^{\N,\NNN}_{2,p,1}$, classical results for the transport equation (see \cite{BCD}, Chap 3) ensure that $h^{2}\in
\widetilde{C}(\R^{+},\widetilde{B}^{\N,\NNN}_{2,p,1})$. And as previously, we have shown that $h^{2}\in
\widetilde{C}(\R^{+},\widetilde{B}^{\N-1,\NNN-1}_{2,p,1})$, it means clearly that $q\in
\widetilde{C}(\R^{+},\widetilde{B}^{\N-1,\NNN}_{2,p,1})$.\\
For getting the continuity result for $u^{2}$, one follows the same ideas.
\subsection*{The proof of the uniqueness}
In the case $\frac{2}{N}\leq\frac{1}{p}+\frac{1}{p_{1}}$, the uniqueness has been established in \cite{BCD}.
\subsection{Proof of the theorems \ref{theo3}}
We follow here exactly the lines
of the proof of the theorem \ref{theo1} except that we take profit that $(\rho^{1},-\mu\n\ln\rho^{1})$ is an exact solution of our system (\ref{10.1}). Indeed we have no coupling by the pressure term in the momentum equation on $u^{2}$, it explains why we do not assume any smallness condition on $q^{1}_{0}$.
\subsection{Proof of the corollary \ref{theo2}}
We remind the reader that if $K_{\mu}$  is the standard heat kernel, that is the solution of the following system:
$$
\begin{cases}
\begin{aligned}
&\p_{t}K_{\mu}-\mu\D K_{\mu}=0,\\
&K_{\mu}(0,\cdot)=\delta_{0}
\end{aligned}
\end{cases}
$$
then:
\begin{equation}
\|D^{\alpha}K_{\mu}(t,\cdot)\|_{L^{p}(\R^{N})}\leq C(\alpha,\mu)t^{-r_{\alpha,p}}
\label{imp4}
\end{equation}
with $r_{\alpha,p}=\N(1-\frac{1}{p})+\frac{\alpha}{2}$. In the sequel we shall note:
$$M_{1}=\|h^{2}\|_{\widetilde{L}^{\infty}(\widetilde{B}^{\N-1,\NNN}_{2,p,1})}\;\;\mbox{and}\;\;M^{2}=\|u^{2}\|_{\widetilde{L}^{\infty}(\widetilde{B}^{\N-1,\NNN-1}_{2,p,1})}\;\;\mbox{such that}\;\;\max(M_{1},M_{2})\leq\e,$$
with $\e$ small enough if we choose the initial data on $(h^{2}_{0},u^{2}_{0})$ small enough.
As $q^{1}$ verifies an heat equation, we have that $q^{1}=K_{\mu}*q^{1}_{0}$.
We obtain then  by using (\ref{imp4}) that:
\begin{equation}
\begin{aligned}
\|(\rho-1)(t,\cdot)\|_{L^{\infty}}&\leq \|q^{1}(e^{h^{2}}-1)(t,\cdot)\|_{L^{\infty}}+\|q^{1}(t,\cdot)\|_{L^{\infty}}+\|(e^{h^{2}}-1)(t,\cdot)\|_{L^{\infty}},\\
&\leq C\big(\frac{\|q^{1}_{0}\|_{L^{1}}}{(1+t)^{\N}}(1+\e)+\e\big) .
\end{aligned}
\label{estim5}
\end{equation}
Similarly we have:
\begin{equation}
\begin{aligned}
\|u(t,\cdot)\|_{B^{-1}_{\infty,\infty}}&\leq \|u^{1}(t,\cdot)\|_{B^{-1}_{\infty,\infty}}+\|u^{2}(t,\cdot)\|_{B^{-1}_{\infty,\infty}},\\
&\leq C\big(\frac{\|q^{1}_{0}\|_{L^{1}}}{(1+t)^{\N+\frac{1}{2}}}+\e)
\end{aligned}
\label{destim5}
\end{equation}
It achieves the proof of the corollary \ref{theo2}.
\section{Appendix}
This section is devoted to the proof of paraproduct estimates which have been used in section $2$ and $3$. They are based on
paradifferential calculus, a tool introduced by J.-M. Bony in \cite{BJM}. The basic idea of paradifferentiel calculus is that
any product of two distributions $u$ and $v$ can be formally decomposed into:
$$uv=T_{u}v+T_{v}u+R(u,v)=T_{u}v+T^{'}_{v}u$$
where the paraproduct operator is defined by $T_{u}v=\sum_{q}S_{q-1}u\D_{q}v$, the remainder operator $R$, by
$R(u,v)=\sum_{q}\D_{q}u(\D_{q-1}v+\D_{q}v+\D_{q+1}v)$ and $T^{'}_{v}u=T_{v}u+R(u,v)$.
\begin{proposition}
\label{hybrid}
Let $p_{1},p_{2},p_{3},p_{4}\in[1,+\infty], (s_{1},s_{2},s_{3},s_{4})\in\R^{4}$ and $(p,q)\in[1,+\infty]^{2}$, we have then the following inequalities:
\begin{itemize}
\item If $\frac{1}{p}\leq\frac{1}{p_{2}}+\frac{1}{\lambda}\leq 1$, $\frac{1}{q}\leq\frac{1}{p_{4}}+\frac{1}{\lambda^{'}}\leq 1$ with $(\lambda,\lambda^{'})\in[1,+\infty]^{2}$ and $p_{1}\leq\lambda^{'}$, $p_{1}\leq\lambda$, $p_{3}\leq\lambda^{'}$ then:
\begin{equation}
\|T_{u}v\|_{\widetilde{B}^{s_{1}+s_{2}+\NN-\frac{N}{p_{1}}-\frac{N}{p_{2}}
,s_{3}+s_{4}+\frac{N}{q}-\frac{N}{p_{3}}-\frac{N}{p_{4}}}_{p,q,r}}\lesssim \|u\|_{\widetilde{B}^{s_{1},s_{3}}_{p_{1},p_{3},1}}\|v\|_{\widetilde{B}^{s_{2},s_{4}}_{p_{2},p_{4},r}},
\label{62}
\end{equation}
if $s_{1}+\frac{N}{\lambda^{'}}\leq\frac{N}{p_{1}}$, $s_{1}+\frac{N}{\lambda}\leq\frac{N}{p_{1}}$ and $s_{3}+\frac{N}{\lambda^{'}}\leq\frac{N}{p_{3}}$.
\item If 
$\frac{1}{q}\leq\frac{1}{p_{3}}+\frac{1}{p_{4}}$ and
$s_{3}+s_{4}+N\inf(0,1-\frac{1}{p_{3}}-\frac{1}{p_{4}})>0$
then
 \begin{equation}
 \sum_{l\geq 4}2^{l(s_{3}+s_{4}+\frac{N}{q}-\frac{N}{p_{3}}-\frac{N}{p_{4}})}\|\D_{l}R(u,v)\|_{L^{q}}\lesssim \|u\|_{\widetilde{B}^{s_{1},s_{3}}_{p_{1},p_{3},1}}\|v\|_{\widetilde{B}^{s_{2},s_{4}}_{p_{2},p_{4},r}}.
\label{63}
\end{equation}
\item If 
$\frac{1}{p}\leq\frac{1}{p_{3}}+\frac{1}{p_{4}}\leq 1$, $\frac{1}{p}\leq\frac{1}{p_{3}}+\frac{1}{p_{2}}\leq 1$, $\frac{1}{p}\leq\frac{1}{p_{1}}+\frac{1}{p_{4}}\leq 1$, $\frac{1}{p}\leq\frac{1}{p_{1}}+\frac{1}{p_{2}}\leq 1$ and
$s_{3}+s_{4}>0$, $s_{3}+s_{2}>0$, $s_{4}+s_{1}>0$, $s_{1}+s_{2}>0$
then
 \begin{equation}
 \sum_{l\leq 4}2^{l(s_{1}+s_{2}+\frac{N}{p}-\frac{N}{p_{1}}-\frac{N}{p_{2}})}\|\D_{l}R(u,v)\|_{L^{p}}\lesssim \|u\|_{\widetilde{B}^{s_{1},\frac{N}{p_{3}}-\frac{N}{p_{1}}+s_{1}}_{p_{1},p_{3},1}}\|v\|_{\widetilde{B}^{s_{2},\frac{N}{p_{4}}-\frac{N}{p_{2}}+s_{2}}_{p_{2},p_{4},r}}.
\label{63}
\end{equation}
with $s_{3}=\frac{N}{p_{3}}-\frac{N}{p_{1}}+s_{1}$ and $s_{4}=\frac{N}{p_{4}}-\frac{N}{p_{2}}+s_{2}$.
\item If $u\in L^{\infty}$, we also have:
\begin{equation}
\|T_{u}v\|_{\widetilde{B}^{s_{1}
,s_{2}}_{p,q,r}}\lesssim \|u\|_{L^{\infty}}\|v\|_{\widetilde{B}^{s_{1},s_{2}}_{p,q,r}},
\label{65}
\end{equation}
and if $\min(s_{1},s_{2})>0$ then:
\begin{equation}
\|R(u,v)\|_{\widetilde{B}^{s_{1}
,s_{2}}_{p,q,r}}\lesssim \|u\|_{L^{\infty}}\|v\|_{\widetilde{B}^{s_{1},s_{2}}_{p,q,r}}.
\label{66}
\end{equation}
\end{itemize}
\end{proposition}
{\bf Proof:} Let us prove (\ref{62}). According to the decomposition of J.-M. Bony \cite{BJM}, we have:
$$uv=T_{u}v+T_{v}u+R(u,v),$$
so for all $l>0$:
$$\D_{l}T_{u}v=\sum_{|l-l^{'}|\leq3}\D_{l}(S_{l^{'}-1}u\D_{l^{'}}v),$$
For $\alpha,\beta\in\R$, let us define the following characteristic function on $\mathbb{Z}$
$$
\begin{aligned}
\va^{\alpha,\beta}=\alpha\;\;\;\mbox{if}\;\;r\leq0,\\
\va^{\alpha,\beta}=\beta\;\;\;\mbox{if}\;\;r\geq1.
\end{aligned}
$$
if  $\frac{1}{p}\leq\frac{1}{p_{2}}+\frac{1}{\lambda}\leq 1$ and $\frac{1}{q}\leq\frac{1}{p_{4}}+\frac{1}{\lambda^{'}}\leq 1$ then
$$\|\D_{l}T_{u}v\|_{L^{\va^{p,q}(l)}}\lesssim 2^{lN\va^{\frac{1}{p_{2}}+\frac{1}{\lambda}-\frac{1}{p},\frac{1}{p_{4}}+\frac{1}{\lambda^{'}}-\frac{1}{q}}(l)}
\sum_{|l-l^{'}|\leq3}\|S_{l^{'}-1}u\|_{L^{\va^{\lambda,\lambda^{'}}(l^{'})}}\|
\D_{l^{'}}v\|_{L^{\va^{p_{2},p_{4}}(l^{'})}}.$$
We have by Berstein inequalities and as $p_{1}\leq\lambda^{'}$, $p_{3}\leq\lambda^{'}$, $p_{1}\leq\lambda$ and
$s_{1}+\frac{N}{\lambda}\leq \frac{N}{p_{1}}$, $s_{1}+\frac{N}{\lambda^{'}}\leq \frac{N}{p_{1}}$, $s_{3}+\frac{N}{\lambda^{'}}\leq \frac{N}{p_{3}}$:
$$
\begin{aligned}
\|S_{l^{'}-1}u\|_{L^{\va^{\lambda,\lambda^{'}}(l^{'})}}&\lesssim\sum_{k\leq l^{'}-2}2^{k(\va^{\frac{N}{p_{1}},\frac{N}{p_{3}}}(k)-\va^{\frac{N}{\lambda},\frac{N}{\lambda^{'}}}(l^{'}))}
\|\D_{k}u\|_{L^{\va^{p_{1},p_{3}}(k)}},\\
&\lesssim\sum_{k\leq l^{'}-2}2^{k(\va^{\frac{N}{p_{1}}-s_{1},\frac{N}{p_{3}}-s_{3}}(k)-\va^{\frac{N}{\lambda},\frac{N}{\lambda^{'}}}(l^{'}))}2^{k\va^{s_{1},s_{3}}(k)}
\|\D_{k}u\|_{L^{\va^{p_{1},p_{3}}(k)}},\\
&\lesssim 2^{l^{'}(\va^{\frac{N}{p_{1}}-s_{1},\frac{N}{p_{3}}-s_{3}}(l^{'})-\va^{\frac{N}{\lambda},\frac{N}{\lambda^{'}}}(l^{'}))}
\|u\|_{\widetilde{B}^{s_{1},s_{3}}_{p_{1},p_{3},1}}.
\end{aligned}
$$
Since $\|\D_{l^{'}}v\|_{L^{\va^{p_{2},p_{4}}(l^{'})}}=c_{l^{'}}2^{-l^{'}(\va^{s_{2},s_{4}}(l^{'}))}
\|v\|_{\widetilde{B}^{s_{2},s_{4}}_{p_{2},p_{4},1}}$ with $\sum_{l^{'}\in\mathbb{Z}}c_{l^{'}}\leq1$
we finally gather as $l>0$:
$$
\begin{aligned}
\|\D_{l}T_{u}v\|_{L^{q}}\lesssim c_{l}2^{l \va^{\frac{N}{p_{1}}+\frac{N}{p_{2}}-\frac{N}{p}-s_{1}-s_{2},
\frac{N}{p_{2}}+\frac{N}{p_{4}}-\frac{N}{q}-s_{3}-s_{4}}(l)}
\|u\|_{\widetilde{B}^{s_{1},s_{3}}_{p_{1},p_{3},1}}\|v\|_{\widetilde{B}^{s_{2},s_{4}}_{p_{2},p_{4},1}}.
\end{aligned}
$$
And we obtain (\ref{62}).
\\
\\
Straightforward modification give (\ref{65}). In this case as $\|S_{k-1}u\|_{L^{\infty}}\leq \|u\|_{L^{\infty}}$ we have:
$$\|\D_{l}T_{u}v\|_{L^{\va^{p,q}(l)}}\lesssim \sum_{|l-l^{'}|\leq3}\|u\|_{L^{\infty}}\|
\D_{l^{'}}v\|_{L^{\va^{p_{2},p_{4}}(l^{'})}}.$$
Next we have:
$$2^{l\va^{p_{2},p_{4}}(l)}\|\D_{l}T_{u}v\|_{L^{\va^{p,q}(l)}}\lesssim \|u\|_{L^{\infty}} \sum_{|l-l^{'}|\leq3}2^{l\va^{p_{2},p_{4}}(l)-l^{'}\va^{p_{2},p_{4}}(l^{'}))}2^{\va^{p_{2},p_{4}}(
l^{'})}\|
\D_{l^{'}}v\|_{L^{\va^{p_{2},p_{4}}(l^{'})}}.$$
We conclude by convolution.\\
\\
To prove (\ref{63}), we write:
$$\D_{l}R(u,v)=\sum_{k\geq l-2}\D_{l}(\D_{k}u\widetilde{\D}_{k}v).$$
We consider now the case $l>3$. By Bernstein and H\"older inequalities we obtain when  $\frac{1}{q}\leq \frac{1}{p_{3}}+\frac{1}{p_{4}}\leq 1$:
$$\|\D_{l}R(u,v)\|_{L^{q}}\lesssim2^{Nl(\frac{1}{p_{3}}+\frac{1}{p_{4}}-\frac{1}{q})}\sum_{k\geq l-2}
\|\D_{k}u\|_{L^{p_{3}}}\|\widetilde{\D}_{k}v\|_{L^{p_{4}}}.$$
Next we have:
$$
\begin{aligned}
2^{l(s_{3}+s_{4}+\frac{N}{q}-\frac{N}{p_{3}}-\frac{N}{p_{4}})}\|\D_{l}R(u,v)\|_{L^{q}}&\lesssim
\sum_{k\geq l-2}
2^{(l-k)(s_{3}+s_{4})}2^{ks_{3}}\|\D_{k}u\|_{L^{p_{3}}}2^{ks_{4}}\|\widetilde{\D}_{k}v\|_{L^{p_{4}}},\\
&\lesssim (c_{k})*(d_{k^{'}}),
\end{aligned}
$$
with $c_{k}=1_{[-\infty,2]}(k)2^{k(s_{3}+s_{4})}$ and $d_{k^{'}}=2^{k^{'}s_{3}}\|\D_{k}u\|_{L^{p_{3}}}2^{k^{'}s_{4}}\|\widetilde{\D}_{k}v\|_{L^{p_{4}}}$
. We conclude by Young inequality as $s_{3}+s_{4}>0$.\\
\\
We have to treat now the case when $l<0$. We have then as $\frac{1}{p}\leq\frac{1}{p_{3}}+\frac{1}{p_{4}}\leq 1$ and
$\frac{1}{p}\leq\frac{1}{p_{1}}+\frac{1}{p_{2}}\leq 1$:
$$
\begin{aligned}
&\|\D_{l}R(u,v)\|_{L^{p}}\lesssim2^{Nl(\frac{1}{p_{3}}+\frac{1}{p_{4}}-\frac{1}{p})}\sum_{k\geq 2}
\|\D_{k}u\|_{L^{p_{3}}}\|\widetilde{\D}_{k}v\|_{L^{p_{4}}}\\
&+\sum_{0\leq k\leq 1, |k-k^{'}|\leq 1}
\|\D_{k}u\D_{k^{'}}v\|_{L^{p}}+2^{Nl(\frac{1}{p_{1}}
+\frac{1}{p_{2}}-\frac{1}{p})}\sum_{l-2\leq k\leq -1}
\|\D_{k}u\|_{L^{p_{1}}}\|\widetilde{\D}_{k}v\|_{L^{p_{2}}}.
\end{aligned}
$$
And by convolution on the middle frequencies:
$$
\begin{aligned}
&\|\D_{l}R(u,v)\|_{L^{p}}\lesssim2^{Nl(\frac{1}{p_{3}}+\frac{1}{p_{4}}-\frac{1}{p})}\sum_{k\geq 2}
\|\D_{k}u\|_{L^{p_{3}}}\|\widetilde{\D}_{k}v\|_{L^{p_{4}}}+(2^{l(\frac{N}{p_{3}}
+\frac{N}{p_{2}}-\frac{N}{p}-s_{3}-s_{2})}c_{l}\\
&\hspace{2cm}+2^{l(\frac{N}{p_{1}}
+\frac{N}{p_{4}}-\frac{N}{p}-s_{1}-s_{4})}
c_{l}+
2^{Nl(\frac{1}{p_{1}}
+\frac{1}{p_{2}}-\frac{1}{p})}\sum_{l-2\leq k\leq -1}
\|\D_{k}u\|_{L^{p_{1}}}\|\widetilde{\D}_{k}v\|_{L^{p_{2}}},
\end{aligned}
$$
with $c_{l}\in l^{1}(\mathbb{Z})$.
Next by convolution we obtain:
$$
\begin{aligned}
&\|\D_{l}R(u,v)\|_{L^{p}}\lesssim c_{l}(2^{l(\frac{N}{p_{3}}+\frac{N}{p_{4}}-\frac{N}{p}-s_{3}-s_{4})}+2^{l(\frac{N}{p_{3}}
+\frac{N}{p_{2}}-\frac{N}{p}-s_{3}-s_{2})}+2^{l(\frac{N}{p_{1}}
+\frac{N}{p_{4}}-\frac{N}{p}-s_{1}-s_{4})}\\
&\hspace{6cm}+
2^{l(\frac{N}{p_{1}}+\frac{N}{p_{2}}-\frac{N}{p}-s_{1}-s_{2})})
\|u\|_{\widetilde{B}^{s_{1},s_{3}}_{p_{1},p_{3},1}}\|v\|_{\widetilde{B}^{s_{2},s_{4}}_{p_{2},p_{4},r}}.
\end{aligned}
$$
And we can conclude.
\\
We want prove now the inequality (\ref{65}). We have then:
$$
\begin{aligned}
2^{l\va^{s_{1},s_{2}}(l)}\|\D_{l}R(u,v)\|_{L^{p}}&\lesssim\sum_{k\geq l-2}2^{(l-k)\va^{s_{1},s_{2}}(l)}
2^{k\va^{s_{1},s_{2}}(l)}\|\D_{k}u\|_{L^{\infty}}\|\widetilde{\D}_{k}v\|_{L^{\va^{p,q}(k)}},
\end{aligned}
$$
And we conclude by Young inequality.
{\hfill $\Box$}

\end{document}